\documentclass[graybox]{svmult}

\usepackage{type1cm}        
\usepackage{makeidx}         
\usepackage{graphicx}        
\usepackage{multicol}        
\usepackage[bottom]{footmisc}

\usepackage{newtxtext}       %
\usepackage{newtxmath}       
\usepackage{algorithm,algorithmic}

\usepackage{wrapfig}

\makeindex             

\usepackage{booktabs}

\usepackage{subcaption}
\usepackage{changepage} 

\usepackage{siunitx}

\usepackage{lineno}

\usepackage{color}
\usepackage{tikz}
\usepackage{pgfplots}
\pgfplotsset{compat=1.10}
\usetikzlibrary{patterns}
\usepackage{pgfplotstable}
\usepackage{csvsimple}
\definecolor{excelblue}{RGB}{94,156,211}
\definecolor{excelorange}{RGB}{235,125,60}
\definecolor{excelgray}{RGB}{165,165,165}
\definecolor{excelgreen}{RGB}{114,172,77}

\usepackage{xcolor,colortbl}
\definecolor{Gray}{gray}{0.9}
\newcolumntype{G}{>{\columncolor{Gray}}r}

\definecolor{tablecolor1}{HTML}{F0F0F0}
\definecolor{tablecolor2}{HTML}{FFFFFF}

\definecolor{mycolor0}{HTML}{2B83BA}
\definecolor{mycolor1}{HTML}{D7301F}
\definecolor{mycolor2}{HTML}{FC8D59}
\definecolor{mycolor3}{HTML}{35A34A}
\definecolor{mycolor4}{HTML}{AF5CBC}
\definecolor{mycolor5}{HTML}{FEF0D9}

\newcommand{\eref}[1]{(\ref{#1})}

\newcommand{\Set}[1]{\mathbb{#1}}
\newcommand{\R}[1]{\Set{R}^{#1}}
\newcommand{\ComplexSet}[1]{\Set{C}^{#1}}

\makeatletter
\newcommand{\leqnos}{\tagsleft@true\let\veqno\@@leqno}
\newcommand{\reqnos}{\tagsleft@false\let\veqno\@@eqno}
\reqnos
\makeatother

\newcommand{\pmat}[1]{\begin{pmatrix}#1\end{pmatrix}}

\def\b{{\discr{b}}}

\def\bfd{{\discr{d}}}
\def\bfe{{\discr{e}}}

\def\bfp{{\mathbf p}}
\def\bfq{{\mathbf q}}
\def\r{{\discr{r}}}
\def\bfs{{\discr{s}}}
\def\bfu{{\discr{u}}}
\def\bfv{{\mathbf v}}

\def\bfx{{\discr{x}}}
\def\x{{\discr{x}}}
\def\bfy{{\discr{y}}}
\def\y{{\discr{y}}}
\def\bfz{{\discr{z}}}

\newcommand{\Complex}[1]{\underline{#1}}
\def\Ybc{\Complex{\mathbf{Y}}^{\text{B}}}
\def\vc{\Complex{\mathbf{v}}}
\def\ic{\Complex{\mathbf{I}}}
\def\sc{\Complex{\mathbf{S}}}
\def\sdc{\Complex{\mathbf{S}}^D}
\def\sgc{\Complex{\mathbf{S}}^G}
\def\sfc{\Complex{\mathbf{S}}^f}
\def\stc{\Complex{\mathbf{S}}^t}

\def\ce{\discr{c}_\varepsilon}

\newcommand{\ciX}[1]{\discr{c}_{\mathcal{I}#1}}
\def\ne{N_\varepsilon}
\def\Je{\discr{J}_\varepsilon}
\newcommand{\JeXX}[2]{\discr{J}_{\varepsilon#1}^{#2}}
\def\lamE{\boldsymbol{\lambda}_\varepsilon}
\def\ci{\discr{c}_\mathcal{I}}
\def\ni{N_\mathcal{I}}
\def\Ji{\discr{J}_\mathcal{I}}
\newcommand{\JiXX}[2]{\discr{J}_{\mathcal{I}#1}^{#2}}
\def\lamI{\boldsymbol{\lambda}_\mathcal{I}}
\def\lamS{\boldsymbol{\lambda}_\bfs}
\def\lamX{\boldsymbol{\lambda}_\bfx}
\def\lamA{\boldsymbol{\lambda}_A}

\def\bfS{{\discr{S}}}

\def\O{{\discr{O}}}
\def\Sc{{\discr{S}}}

\def\A{{\discr{A}}}
\def\B{{\discr{B}}}

\def\Bs{{\mathbf{B}^\text{S}}}
\def\E{{\mathbf {E}}}
\def\F{{\mathbf {F}}}
\def\C{{\discr{C}}}
\def\M{{\discr{M}}}

\def\X{{\discr{X}}}

\def\L{{\discr{L}}}

\def\PowerFlow{\discr{c}_\varepsilon}
\def\LineFlow{\discr{c}_\mathcal{I}}

\def\PowerFlowSecurity{{\PowerFlow^c}}
\def\LineFlowSecurity{{\LineFlow^c}}

\newcommand{\pes}[2]{{\mathbf p}^{#1}_{#2}}

\def\psmax{{\bfp}^{\text{S,max}}}
\def\psmin{{\bfp}^{\text{S,min}}}
\def\qsmax{{\bfq}^{\text{S,max}}}
\def\qsmin{{\bfq}^{\text{S,min}}}

\def\pg{{\bfp}^{\text{G}}}
\def\qg{{\bfq}^{\text{G}}}
\def\ps{{\bfp}^{\text{S}}}
\def\qs{{\bfq}^{\text{S}}}
\newcommand{\psci}{{\bfp}^{\text{Sc},i}}
\newcommand{\psdi}{{\bfp}^{\text{Sd},i}}

\newcommand{\pscni}[2]{{\bfp}^{\text{Sc},{#1}}_{#2}}
\newcommand{\psdni}[2]{{\bfp}^{\text{Sd},{#1}}_{#2}}
\newcommand{\qscni}[2]{{\bfq}^{\text{Sc},{#1}}_{#2}}
\newcommand{\qsdni}[2]{{\bfq}^{\text{Sd},{#1}}_{#2}}

\def\scomplex{\Complex{\mathbf S}}

\def\sg{\Complex{\mathbf S}^{\text{G}}}
\def\sS{\Complex{\mathbf S}^{\text{S}}}

\def\btheta{{\boldsymbol \theta}}

\def\bflambda{{\boldsymbol \lambda}}
\def\bepsilon{{\boldsymbol \epsilon}}

\def\0{{\mathbf 0}}

\def\Nb{{N_\text{B}}}
\def\Ng{{N_\text{G}}}
\def\Ns{{N_\text{S}}}
\def\Nl{{N_\text{L}}}

\def\Cgb{C^{\text{G}}}

\newcommand{\Hess}{\discr{H}}
\newcommand \discr[1]{{\boldsymbol{#1}}}   

\newcommand{\BELTISTOS}{{\small BELTISTOS}}
\newcommand{\MATPOWER}{{\small MATPOWER}}
\newcommand{\PARDISO}{{\small PARDISO}}
\newcommand{\IPOPT}{{\small IPOPT}}
\newcommand{\MIPS}{{\small MIPS}}
\newcommand{\KNITRO}{{\small KNITRO}}

\newcommand \Min[1]{{#1}^{\text{min}}}
\newcommand \Max[1]{{#1}^{\text{max}}}

\def\grad{\nabla}

\def\superstar{^{\raise 0.5pt\hbox{$\nthinsp *$}}}

\def\nthinsp{\mskip -2   mu}
\makeatletter
\newcommand*{\transpose}{%
	{\mathpalette\@transpose{}}%
}
\newcommand*{\@transpose}[2]{%
	\raisebox{\depth}{$\m@th#1\intercal$}%
}
\makeatother

\def\minim{\mathop{\hbox{\rm minimize}}}
\def\subject{\text{\rm subject to}}
\def\minimize#1{{\displaystyle\minim_{#1}}}

\begin{document}

\title*{Structure-Exploiting Interior Point Methods}
\author{Juraj Kardo\v s \and Drosos Kourounis \and Olaf Schenk}
\authorrunning{J. Kardo\v s, D. Kourounis, O. Schenk}
\institute{Juraj Kardo\v s \at Institute of Computational Science, Faculty of Informatics, Universit\'a della Svizzera italiana,
         Switzerland, \email{juraj.kardos@usi.ch}
\and Drosos Kourounis \at Institute of Computational Science, Faculty of Informatics, Universit\'a della Svizzera italiana,
         Switzerland, \email{drosos.kourounis@usi.ch}         
\and Olaf Schenk \at Institute of Computational Science, Faculty of Informatics, Universit\'a della Svizzera italiana,
         Switzerland, \email{olaf.schenk@usi.ch}}
%
%
\maketitle

\abstract*{
}
\abstract{
Interior point methods are among the most popular techniques for large-scale nonlinear 
optimization, owing to their intrinsic ability of scaling to arbitrarily large 
problem sizes. Their efficiency has attracted, in recent years, a lot of  attention due to increasing demand for large-scale optimization in industry and engineering. General purpose nonlinear programming solvers implementing interior point methods provide a generic framework that embraces almost every 
optimal control problem provided that the nonlinear functions are smooth and that their first and possibly the second derivatives are also available. However, a large class of industrial and engineering problems possesses a particular structure motivating the development of structure-exploiting interior point methods. We present an interior point framework that exploits the intrinsic structure of large-scale nonlinear optimization problems, with the purpose of accelerating the solution both for single-core or multicore execution and massively parallel high-performance computing infrastructures. 
Since the overall performance of interior point methods relies heavily on scalable sparse linear algebra solvers, particular emphasis is given to the underlying algorithms for the distributed solution of the associated sparse linear systems obtained at each iteration from the linearization of the optimality conditions. The interior point algorithm is implemented in an  object-oriented parallel solver and applied for the solution of large-scale optimal control problems solved on a daily basis for the secure transmission and distribution of electricity in modern power grids.}


\section{Introduction}
\label{sec:intro}

Interior point (IP) methods  have became a successful tool for solving the nonlinearly constrained optimal control problems. Their origin can be traced back to  1984 when Karmarkar \cite{Karmarkar1984} announced a polynomial time linear program that was considerably faster than the most popular simplex method to date. Furthermore, IP methods can also be applied to quadratic and other nonlinear programs, unlike the simplex method which can be applied only to linear programming. 
The main advantages of the IP methods lie in the convenience they offer for handling nonlinear inequality constraints using logarithmic barrier functions, so that a strictly feasible initial point is unnecessary.
Another advantage of IP methods is that they are applicable to large-scale problems and allow for a variety of different direct sparse or iterative solution methods for the underlying linear systems solved at each iteration until convergence. Since different sparse system solvers can be plugged in with ease, large-scale structured problems can be solved by exploiting parallel computing infrastructures. 


An example of successful application of the IP methods is the class of problems known as the optimal power flow (OPF). OPF is a nonlinear, nonconvex, large-scale optimization problem with the objective of minimizing the electricity generation cost while satisfying the physical constraints of the electric grid. The security constrained OPF (SCOPF) is an extension of the OPF problem that additionally ensures the system security with respect to a set of postulated contingencies. The SCOPF has become an essential tool for many transmission system operators for the planning, operational planning, and real time operation of the power system. An  increase  of  the number of considered contingencies requires the introduction of  additional  variables  and  constraints,  which  in  turn  results  in  a significant problem size growth, rendering the solution computationally intractable  for  standard  general  purpose  optimization  tools. The structure of the SCOPF problems is appropriate for the parallel structure exploiting IP methods, where each contingency corresponds to a separate partition on the linear level. 
The nonlinear IP framework leverages the bordered block-diagonal sparse structure specific to these optimal control problems by applying a Schur complement elimination on a  block-per-block basis in order to exploit parallelism intrinsic to sparse block-diagonal structures by distributing the block contributions to the global Schur complement. 
In this way, the solution of the large-scale optimization problems can be approached more efficiently, as demonstrated in~\cite{CosminMPSCOPF}. Similar structures arise also in the multistage stochastic optimal control problems \cite{CosminCurvature, pardiso}, multiperiod OPF problems (MPOPF)~\cite{beltistos}, dynamic simulations of the power grid \cite{flexdyn}, or problems such as natural gas dispatch \cite{CosminCurvature}.

This overview summarizes the algorithmic improvements in the recent years that have significantly advanced IP methods. The focus is on parallel implementations demonstrated on problems arising from the optimal control of the power grid. The presented primal-dual IP method is based on the \IPOPT{} algorithm \cite{ipoptThesis, IPOPT2006}.

\subsection{Notation}
Throughout we adopt the following notation. Scalar values are denoted by  lowercase letters $x$ in normal font, while vector objects are represented by bold lowercase letters $\bfx$. The vector $\bfe$ is a vector of ones with an appropriate dimension. If not specified otherwise, column vectors are assumed. Similarly, scalar functions are represented by a lowercase letter $f$, while vector functions are shown in bold lowercase $\boldsymbol{f}$. Concatenation of column vectors $(\bfx_1^\transpose,\bfx_2^\transpose, \ldots )^\transpose$ will be denoted by $(\bfx_1,\bfx_2, \ldots )$. The elementwise product of two vectors $\x,\y$ will be denoted by $\bfx \bfy$, while $\bfx^\transpose \bfy$ stands for the inner product of the two vectors. Matrices are represented by uppercase letters; for general (sparse) matrices we use bold fonts $\X$ while we will use normal font to distinguish diagonal matrices $X$. Sets will be represented by a calligraphic font $\mathcal{X}$ or uppercase Greek letters.

\section{IP Algorithm\label{chapter:IPM}}
\begin{definition}
A general nonlinear programming (NLP) problem is formulated as a minimization problem
\begin{subequations}
	\label{eq:NLP}
	\begin{align}
	\minimize{\bfx} \; & f(\bfx)  \label{eq:NLPobj}\\
	\text{subject to} \; &  \ce(\bfx) = \0,  \label{eq:NLPeq}\\
	&  \ci(\bfx) \geq \0,  \label{eq:NLPineq}\\
	& \bfx \geq \0, \label{eq:NLPbounds}
	\end{align}
\end{subequations}
where $\bfx \in \mathbb{R}^{N_x}$, the objective function $f$ is a mapping $f: \mathbb{R}^{N_x} \rightarrow \mathbb{R}$, the constraints $\ce: \mathbb{R}^{N_x} \rightarrow \mathbb{R}^{\ne}$ and $\ci: \mathbb{R}^{N_x} \rightarrow \mathbb{R}^{\ni}$  are assumed to be sufficiently smooth, with continuous second order derivatives,  and ${N_x} > \ne, \ni$, where $\ne, \ni$ are the number of equality and inequality constraints, respectively.
\end{definition}
%


\begin{definition}
The \emph{feasible set} $\Omega$ is a set of points $\bfx$ that satisfy the constraints of the NLP problem \eqref{eq:NLP}; that is
\begin{equation}
    \Omega = \{\bfx \in \mathbb{R}^{N_x} \ | \ \ce(\bfx) = \0,\ \ci(\bfx) \geq \0,\ \bfx \geq \0 \}.
\end{equation}
\end{definition}

\begin{definition}
The \emph{active set} at any feasible point $\bfx$ is a set of inequality constraints indices, for which the equality constraint holds; that is, $\mathcal{A}(\bfx) = \{i \ |\ \ci^i(\bfx) = \0 \}$.
\end{definition}

\begin{definition}
Given the solution of the NLP problem $\bfx*$ and the active set $\mathcal{A}(\bfx^*)$, the \emph{linear independence constraint qualification} (LICQ) holds if the set of active constraint gradients $\{ \nabla \ce^i(\bfx^*),\ i=1\ldots \ne; \quad \nabla \ci^j(\bfx^*),\  j \in \mathcal{A}(\bfx^*)\}$, is linearly independent.
\end{definition}
The NLP problem \eqref{eq:NLP} can be transformed into the equivalent problem formulation where the inequality constraints are converted to equality constraints by introducing the slack variables $\bfs \in \mathbb{R}^{\ni}$ with additional nonnegativity bounds $\bfs \geq \0$. The NLP problem can be written as
\begin{subequations}
	\label{eq:NLPslack}
	\begin{align}
	\minimize{\bfx} \; & f(\bfx)  \label{eq:NLPslackobj}\\
	\text{subject to} \; &  \ce(\bfx) = \0,  \label{eq:NLPslackeq}\\
	&  \ci(\bfx) - \bfs = \0,  \label{eq:NLPslackineq}\\
	& (\bfx, \bfs) \geq \0. \label{eq:NLPslackbounds}
	\end{align}
\end{subequations}


\begin{definition}
The \emph{Lagrangian} for the NLP problem \eqref{eq:NLPslack} is defined as 
\begin{equation}
    \label{eq:lagrangianNLP}
    \mathcal{L}(\bfx, \bfs, \lamE, \lamI, \bflambda_x, \bflambda_s) = f(\bfx) + \lamE^\transpose \ce(\bfx)  + \lamI^\transpose (\ci(\bfx) - \bfs) - \bflambda_x^\transpose \bfx - \bflambda_s^\transpose \bfs.
\end{equation}
\end{definition}
The vectors $\lamE, \lamI, \bflambda_x$, and $\bflambda_s$ are the Lagrange multipliers associated with the equality, original inequality, and the bound constraints on the primal and slack variables. This allows us to state the Karush--Kush--Tucker (KKT) first-order necessary conditions for the NLP problem \eqref{eq:NLPslack} which characterize the solution.
%
%
\begin{theorem}
\label{theorem:KKT} Suppose that $\bfx^*$ is a local solution of the NLP problem \eqref{eq:NLPslack} and that the LICQ holds at $\bfx^*$. Then there exist Lagrange multiplier vectors $\lamE^* \in \mathbb{R}^{\ne}, \lamI^* \in \mathbb{R}^{\ni}, \bflambda_x^* \in \mathbb{R}^n$ and $\bflambda_s^* \in \mathbb{R}^{\ni}$, $(\lamX^*, \lamS^*) \geq \0$, such that the following conditions are satisfied at $(\bfx^*, \bfs^*, \lamE^*, \lamI^*, \bflambda_x^*, \bflambda_s^*)$:
\begin{subequations}
\label{eq:KKT}
\begin{align}
    \grad_x f(\bfx^*) + \grad_x \ce(\bfx^*)^\transpose \lamE^{*} + \grad_x \ci(\bfx^*)^\transpose \lamI^{*} - \bflambda_x^{*} &= \0, \label{eq:KKTdual_x}\\
    - \lamI^{*} - \bflambda_s^{*} &= \0, \label{eq:KKTdual_s}\\
    \ce (\bfx^*) &= \0, \label{eq:KKTprimalG}\\
    \ci (\bfx^*) - \bfs^* &= \0, \label{eq:KKTprimalH}\\
    \lamX^* \bfx^* &= \0,  \label{eq:KKTcomplementX} \\
    \lamS^* \bfs^* &= \0, \label{eq:KKTcomplementH}\\
    (\bfx^*, \bfs^*) &\geq \0.
\end{align}
\end{subequations}
\end{theorem}
The conditions \eqref{eq:KKTdual_x} and \eqref{eq:KKTdual_s} are referred to as dual feasibility, \eqref{eq:KKTprimalG}, \eqref{eq:KKTprimalH} as primal feasibility, and \eqref{eq:KKTcomplementH}, \eqref{eq:KKTcomplementX} as complementarity conditions. The point $\bfx^*$ satisfying the KKT conditions is called a stationary, or critical, point. In order to ensure that any stationary point $\bfx^*$ is indeed an optimal (local) solution of the NLP problem \eqref{eq:NLPslack}, the second-order sufficient conditions are needed. 
\begin{theorem}
\label{theorem:SOC}
Let $\bfx^*$ be a point at which LICQ holds, the KKT conditions are satisfied, and strict complementarity holds for the active inequality constraints. Then, the point $\bfx^*$ satisfies the second-order sufficient conditions for the NLP problem \eqref{eq:NLPslack} if the Hessian of the Lagrangian $\nabla^2_{xx} \mathcal{L}(\bfx^*, \bfs^*, \lamE^*, \lamI^*, \bflambda_x^*, \bflambda_s^*)$  projected 
onto the null space of the constraint Jacobian is positive definite.
\end{theorem}
In practice, the second-order conditions are guaranteed by monitoring the inertia of the iteration matrix, which is further elaborated in section \ref{sec:inertia}.
Proofs of Theorems \ref{theorem:KKT} and \ref{theorem:SOC} can be found in classic optimization textbooks, e.g., \cite{NocedalBookOptim,WrightPDIPM}.
If the active set at the solution of the NLP problem was known, we could apply a Newton-class method directly to the linearization of the KKT conditions. However, the identification of the active set is known to be an NP-hard combinatorial problem for which, in the worst case, the computation time increases exponentially with the size of the problem. Therefore, many solution strategies adopt an IP approach, introducing a barrier subproblem where the nonnegativity bounds on the variables and slacks $(\bfx, \bfs) \geq \0$ are handled by the standard logarithmic barrier function, which is, in fact, a penalty term penalizing the iterates that approach the boundary of the feasible region. 
\begin{definition}
The \emph{barrier subproblem} (BSP) reads:
\begin{subequations}
	\label{eq:Barrier}
	\begin{align}
	\minimize{\bfx, \bfs} \; & f(\bfx) - \mu \sum_{i=1}^{n} \log (x_i) - \mu \sum_{i=1}^{\ni} \log (s_i)  \label{eq:Barrier_obj}\\
	\text{subject to} \; &  \ce(\bfx) = \0,  \label{eq:Barrier_eq}\\
	&  \ci(\bfx) - \bfs = \0.  \label{eq:Barrier_ineq}
	\end{align}
\end{subequations}
\end{definition}
Under certain conditions the solution $\bfx^*$ of the BSP \eqref{eq:Barrier} converges to the solution of the original NLP problem \eqref{eq:NLP} as  $\mu_j \downarrow 0$. Consequently, a strategy to solve the original NLP problem is to solve a sequence of the BSPs decreasing the barrier parameter $\mu_j$. The solution of each iterate is not relevant for the solution of the original problem, so it can be relaxed to a certain accuracy and such an approximate solution is used as a starting point for the next BSP. The strategy for updating the $\mu$ parameter and thus switching to the next BSP is discussed later in section \ref{sec:Barrier}.

The solutions of the barrier problem \eqref{eq:Barrier} are critical points of the Lagrangian function
\begin{align}
    \label{eq:lagrangian}
  \mathcal{L}(\bfx, \bfs, \lamE, \lamI) =& f(\bfx) -  \mu_j \sum_{i=1}^{N_x} \log(x_i) - \mu_j \sum_{i=1}^{\ni} \log(s_i)\\
  & + \lamE^\transpose \ce(\bfx) + \lamI^\transpose (\ci(\bfx) - \bfs) \nonumber.
\end{align}
Formulating and solving the optimality conditions of \eqref{eq:lagrangian} directly would lead to singularities, since the derivatives of the barrier terms involve the fractions $\frac{\mu}{x_i}$ and $\frac{\mu}{s_i}$, which are not defined at the solution $\bfx^*, \bfs^*$ of the NLP problem \eqref{eq:NLP} when active bounds $x^*_i = 0$ or $s^*_i = 0$ are attained. Primal-dual IP methods~\cite{Conn2000, Conn2001} define the dual variables $\bfz$ and $\bfy$ as
\begin{subequations}
\label{eq:dualVar}
\begin{align}
z_i &= \frac{\mu}{x_i},\ i=1,2,\ldots,N_x, \label{eq:dualVar_x} \\
y_i &= \frac{\mu}{s_i},\ i=1,2,\ldots,\ni. \label{eq:dualVar_slack}
\end{align}
\end{subequations}




From the definition of the dual variables it follows that $z_i = \frac{\mu}{x_i} > 0$; therefore, $z_i x_i = \mu \ \forall i = 1,\ldots,N_x$. Similarly, $y_i s_i = \mu, y_i > 0 \ \forall i = 1,\ldots,\ni$.  The optimality conditions of the BSP \eqref{eq:Barrier}, considering also the dual variables \eqref{eq:dualVar}, are written 
\begin{subequations}
\label{eq:BSPoc}
\begin{align}
    \grad_x f(\bfx^*) + \grad_x \ce(\bfx^*)^\transpose \lamE^{*} + \grad_x \ci(\bfx^*)^\transpose \lamI^{*} -\bfz^* &= \0, \label{eq:BSPoc_dual_x} \\
    - \lamI^* - \bfy^* &= \0, \label{eq:BSPoc_dual_s}\\
    \ce (\bfx^*) &= \0, \label{eq:BSPoc_primalG}\\
    \ci (\bfx^*) -\bfs^* &= \0, \label{eq:BSPoc_primalH}\\
    \bfz^{*} \bfx^* &= \mu \bfe,  \label{eq:BSPoc_complementX} \\
    \bfy^{*} \bfs^* &= \mu \bfe, \label{eq:BSPoc_complementH} \\
    (\bfx^*, \bfs^*) & \geq \0.
\end{align}
\end{subequations}
Note that the dual variables $\bfz, \bfy$ correspond to the Lagrange multipliers $\bflambda_x$ and $\bflambda_s$ for the bound constraints. The KKT conditions of the BSP \eqref{eq:BSPoc} are equivalent to the perturbed conditions \eqref{eq:KKT} of the original NLP problem \eqref{eq:NLPslack}, except for the strict positivity of the dual variables $(\bfz, \bfy) > 0$.
The primal-dual equations then become
\begin{subequations}
	\label{eq:OC}
	 \begin{alignat}{2}
	\boldsymbol{l}_a := \quad & \nabla_x f(\bfx) + \Je^\transpose \lamE + \Ji^\transpose \lamI - \bfz \  &&= \  \0, \\
	\boldsymbol{l}_b := \quad & -\lamI - \bfy \  &&= \  \0 \label{eq:OCsecond}, \\
	\boldsymbol{l}_c := \quad & \ce(\bfx) \  &&= \  \0, \\
	\boldsymbol{l}_d := \quad & \ci(\bfx) - \bfs \  &&= \  \0, \\
	\boldsymbol{l}_e :=  \quad& Z\bfx - \mu \bfe \  &&=  \ \0, \label{eq:OCcomplx}\\
	\boldsymbol{l}_f :=  \quad& Y\bfs - \mu \bfe \  &&= \  \0, \label{eq:OCcompls}
	\end{alignat}
\end{subequations}
where the Jacobian of constraints is written as $\Je = \grad_x \ce(\bfx)$ and $\Ji = \grad_x \ci(\bfx)$. The diagonal matrices $X, S, Z, Y$ are defined as $X = \text{diag}(\bfx)$, $S = \text{diag}(\bfs)$, $Z = \text{diag}(\bfz)$, and $Y = \text{diag}(\bfy)$.

Linearizing the primal-dual equations and solving them by applying Newton's method starting from an arbitrary value of the barrier parameter $\mu$ may result in slow convergence or poor conditioning of the associated KKT systems. Following the central path ensures that certain favorable conditions for the KKT systems and primal-dual variables are satisfied and descent directions can be obtained with reasonable accuracy.

\begin{definition}
The \emph{central path} $\mathcal{C}$ is an arc of strictly feasible points of the BSP problem~\eqref{eq:Barrier}, $\mathcal{C} = \{ (\bfx^{\mu}, \bfs^{\mu}, \lamE^{\mu}, \lamI^{\mu}, \bfz^{\mu}, \bfy^{\mu})\ |\ \mu > 0\}$, such that $(\bfx^{\mu}, \bfs^{\mu}, \lamE^{\mu}, \lamI^{\mu}, \bfz^{\mu}, \bfy^{\mu})$ is a solution of the BSP problem for every value of $\mu > 0$. Points on the central path are characterized by the first-order KKT conditions \eqref{eq:OC}.
\end{definition}
\begin{definition}
 The\emph{ duality measure} $\tau$ is an average pairwise complementarity value $x_i z_i$ and $s_i y_i$,
\begin{equation}
    \label{eq:dualMeasure}
    \tau = \frac{\bfx^\transpose \bfz + \bfs^\transpose \bfy}{N_x+\ni}.
\end{equation}
\end{definition}
The barrier parameter $\mu$ is usually chosen proportionally to the duality measure and the centering parameter $\sigma \in [0,1]$, such that $\mu = \tau \sigma$.
By choosing $\sigma = 1$ the algorithm moves toward the central path $\mathcal{C}$. Such a step is biased toward the interior of the feasible region defined by the constraints $(\bfz, \bfx) > \0$, $(\bfy, \bfs) > \0$. 
At the other extreme, the value $\sigma = 0$ results in the standard Newton step aiming to satisy the KKT conditions \eqref{eq:KKT}. Many algorithms use intermediate values of $\sigma$ from the open interval $(0,1)$ to trade off between the two objectives of reducing duality measure and improving centrality. A strategy for selecting the centering parameter is discussed later in sections \ref{sec:adaptive1} and \ref{sec:adaptive2}.




\begin{remark} \label{remark:boundsVars}
The treatment  for general box constraints $\bfx^{\text{min}} \leq \bfx \leq \bfx^{\text{max}}$  and general upper and lower bounds on the nonlinear constraints  $\ci^{\text{min}} \leq \bfs \leq \ci^{\text{max}}$ requires the addition of modified logarithmic barrier terms 
\begin{align}
    \mathcal{B}(\bfx, \bfx^{\text{min}}, \bfx^{\text{max}}) &= -\mu_j \sum_{i=1}^{N_x} \log(x_i - x_i^{\text{min}}) - \mu_j \sum_{i=1}^{N_x} \log( x_i^{\text{max}} - x_i),
\\
    \mathcal{B}(\bfs, \ci^{\text{min}}, \ci^{\text{max}}) &= -\mu_j \sum_{i=1}^{\ni} \log(s_i - \ciX{i}^{\text{min}}) - \mu_j \sum_{i=1}^{\ni} \log( \ciX{i}^{\text{max}} - s_i).
\end{align}
The dual variables for $i=1,2,\ldots,N_x$ are defined by
\begin{align}
z_i^L &= \frac{\mu}{x_i - x_i^{\text{min}}}, &
z_i^U &= \frac{\mu}{x_i^{\text{max}} - x_i}, &
\label{eq:dualVarFULLz}
\end{align}
while for the constraints the dual variables are defined by
\begin{align}
y_i^L &= \frac{\mu}{s_i - \ciX{i}^{\text{min}}}, & 
y_i^U &= \frac{\mu}{\ciX{i}^{\text{max}} - s_i}. &
\label{eq:dualVarFULLy}
\end{align}
\end{remark}

\subsection{Search Direction Computation  \label{sec:SearchDir}}
Since the solution of the barrier problem \eqref{eq:Barrier} satisfies the perturbed KKT conditions \eqref{eq:OC}, Newton's method may be applied to solve the system of nonlinear equations. The search direction $(\Delta \bfx^k, \Delta \bfs^k, \Delta \lamE^k,\Delta \lamI^k, \Delta \bfz^k, \Delta \bfy^k)$ at the $k$th iteration can be obtained from the linearization of \eqref{eq:OC} at the iterate $( \bfx^k, \bfs^k, \lamE^k, \lamI^k, \bfz^k, \bfy^k)$, resulting in a system of linear equations
\begin{equation}
\label{eq:SCOPFhessian}
\begin{bmatrix}
\boldsymbol{H} & \0 & \Je^\transpose & \Ji^\transpose & -I & \0 \\
\0 & \0 & \0 & -I & \0 & -I \\
\Je & \0 & \0 & \0 & \0 & \0 \\
\Ji & -I & \0 & \0 & \0 & \0\\
Z & \0 & \0 & \0 & X & \0 \\
\0 & Y & \0 & \0 & \0 & S
\end{bmatrix}^k
\begin{bmatrix}
\Delta \bfx\\
\Delta \bfs\\
\Delta \lamE\\
\Delta \lamI\\
\Delta \bfz\\
\Delta \bfy
\end{bmatrix}^k
=
-\begin{bmatrix}
\boldsymbol{l}_a \\
\boldsymbol{l}_b \\
\boldsymbol{l}_c \\
\boldsymbol{l}_d \\
\boldsymbol{l}_e \\
\boldsymbol{l}_f
\end{bmatrix}^k,
\end{equation}
where  $\boldsymbol{H} = \nabla^2_{xx} \mathcal{L}$. The system \eqref{eq:SCOPFhessian} is clearly unsymmetric. A symmetric system can be obtained after eliminating the last two block rows:
\begin{equation}
\label{eq:SCOPFhessianSymmetric}
\begin{bmatrix}
\boldsymbol{\tilde{H}} & \0 & \Je^\transpose & \Ji^\transpose \\
\0 & L_s & \0 & -I \\
\Je & \0 & \0 & \0 & \\
\Ji & -I & \0 & \0 & 
\end{bmatrix}^k
\begin{bmatrix}
\Delta \bfx\\
\Delta \bfs\\
\Delta \lamE\\
\Delta \lamI
\end{bmatrix}^k
=
- \begin{bmatrix}
\boldsymbol{l}_a + {X}^{-1} \boldsymbol{l}_e\\
\boldsymbol{l}_b + {S}^{-1} \boldsymbol{l}_f\\
\boldsymbol{l}_c \\
\boldsymbol{l}_d
\end{bmatrix}^k,
\end{equation}
where $\boldsymbol{\tilde{H}} = \boldsymbol{H} + {X}^{-1}{Z}$ and $L_s = S^{-1}Y$. The directions $\Delta \bfz^k$ and $\Delta \bfy^k$ can then be recovered from the equations
\begin{align}
\label{eq:SCOPFhessianSymmetricRecovery}
    \Delta \bfz^k = -{X}^{-1} (\boldsymbol{l}_e + {Z}\Delta \bfx^k), \\
    \Delta \bfy^k = -{S}^{-1} (\boldsymbol{l}_f + {Y}\Delta \bfs^k).
\end{align}

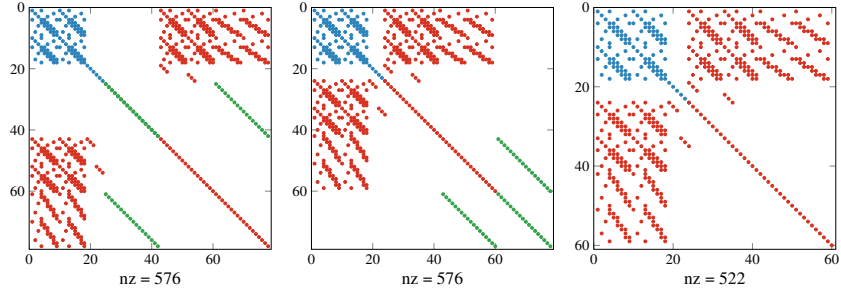
\begin{figure}[t]
\begin{center}
%
%
\begin{tikzpicture}[scale=.55]

\begin{axis}[%
width=0.5\textwidth,
height=0.5\textwidth,
at={(1.648in,0.642in)},
scale only axis,
xmin=0,
xmax=79,
xlabel style={font=\color{white!15!black}},
xlabel={nz = 576},
y dir=reverse,
ymin=0,
ymax=79,
axis background/.style={fill=white},
legend style={legend cell align=left, align=left, draw=white!15!black},
every x tick label/.append style={font=\normalsize},
every y tick label/.append style={font=\normalsize},
label style={font=\large}
]
\addplot [color=black, draw=none, mark size=1pt, mark=*, mark options={solid, mycolor0}]
  table[row sep=crcr]{%
  	    1 1\\
  	1 4\\
  	1 10\\
  	1 13\\
  	2 2\\
  	2 8\\
  	2 11\\
  	2 17\\
  	3 3\\
  	3 6\\
  	3 12\\
  	3 15\\
  	4 1\\
  	4 4\\
  	4 5\\
  	4 9\\
  	4 10\\
  	4 13\\
  	4 14\\
  	4 18\\
  	5 4\\
  	5 5\\
  	5 6\\
  	5 13\\
  	5 14\\
  	5 15\\
  	6 3\\
  	6 5\\
  	6 6\\
  	6 7\\
  	6 12\\
  	6 14\\
  	6 15\\
  	6 16\\
  	7 6\\
  	7 7\\
  	7 8\\
  	7 15\\
  	7 16\\
  	7 17\\
  	8 2\\
  	8 7\\
  	8 8\\
  	8 9\\
  	8 11\\
  	8 16\\
  	8 17\\
  	8 18\\
  	9 4\\
  	9 8\\
  	9 9\\
  	9 13\\
  	9 17\\
  	9 18\\
  	10 1\\
  	10 4\\
  	10 10\\
  	10 13\\
  	11 2\\
  	11 8\\
  	11 11\\
  	11 17\\
  	12 3\\
  	12 6\\
  	12 12\\
  	12 15\\
  	13 1\\
  	13 4\\
  	13 5\\
  	13 9\\
  	13 10\\
  	13 13\\
  	13 14\\
  	13 18\\
  	14 4\\
  	14 5\\
  	14 6\\
  	14 13\\
  	14 14\\
  	14 15\\
  	15 3\\
  	15 5\\
  	15 6\\
  	15 7\\
  	15 12\\
  	15 14\\
  	15 15\\
  	15 16\\
  	16 6\\
  	16 7\\
  	16 8\\
  	16 15\\
  	16 16\\
  	16 17\\
  	17 2\\
  	17 7\\
  	17 8\\
  	17 9\\
  	17 11\\
  	17 16\\
  	17 17\\
  	17 18\\
  	18 4\\
  	18 8\\
  	18 9\\
  	18 13\\
  	18 17\\
  	18 18\\
  	19 19\\
  	20 20\\
  	21 21\\
  	22 22\\
  	23 23\\
  	24 24\\
};

\addplot [color=black, draw=none, mark size=1pt, mark=*, mark options={solid, mycolor1}]
table[row sep=crcr]{%
	    43 1\\
	43 4\\
	43 10\\
	43 13\\
	43 19\\
	44 2\\
	44 8\\
	44 11\\
	44 17\\
	44 20\\
	45 3\\
	45 6\\
	45 12\\
	45 15\\
	45 21\\
	46 1\\
	46 4\\
	46 5\\
	46 9\\
	46 10\\
	46 13\\
	46 14\\
	46 18\\
	47 4\\
	47 5\\
	47 6\\
	47 13\\
	47 14\\
	47 15\\
	48 3\\
	48 5\\
	48 6\\
	48 7\\
	48 12\\
	48 14\\
	48 15\\
	48 16\\
	49 6\\
	49 7\\
	49 8\\
	49 15\\
	49 16\\
	49 17\\
	50 2\\
	50 7\\
	50 8\\
	50 9\\
	50 11\\
	50 16\\
	50 17\\
	50 18\\
	51 4\\
	51 8\\
	51 9\\
	51 13\\
	51 17\\
	51 18\\
	52 1\\
	52 4\\
	52 10\\
	52 13\\
	52 22\\
	53 2\\
	53 8\\
	53 11\\
	53 17\\
	53 23\\
	54 3\\
	54 6\\
	54 12\\
	54 15\\
	54 24\\
	55 1\\
	55 4\\
	55 5\\
	55 9\\
	55 10\\
	55 13\\
	55 14\\
	55 18\\
	56 4\\
	56 5\\
	56 6\\
	56 13\\
	56 14\\
	56 15\\
	57 3\\
	57 5\\
	57 6\\
	57 7\\
	57 12\\
	57 14\\
	57 15\\
	57 16\\
	58 6\\
	58 7\\
	58 8\\
	58 15\\
	58 16\\
	58 17\\
	59 2\\
	59 7\\
	59 8\\
	59 9\\
	59 11\\
	59 16\\
	59 17\\
	59 18\\
	60 4\\
	60 8\\
	60 9\\
	60 13\\
	60 17\\
	60 18\\
	61 1\\
	61 4\\
	61 10\\
	61 13\\
	62 4\\
	62 5\\
	62 13\\
	62 14\\
	63 5\\
	63 6\\
	63 14\\
	63 15\\
	64 3\\
	64 6\\
	64 12\\
	64 15\\
	65 6\\
	65 7\\
	65 15\\
	65 16\\
	66 7\\
	66 8\\
	66 16\\
	66 17\\
	67 2\\
	67 8\\
	67 11\\
	67 17\\
	68 8\\
	68 9\\
	68 17\\
	68 18\\
	69 4\\
	69 9\\
	69 13\\
	69 18\\
	70 1\\
	70 4\\
	70 10\\
	70 13\\
	71 4\\
	71 5\\
	71 13\\
	71 14\\
	72 5\\
	72 6\\
	72 14\\
	72 15\\
	73 3\\
	73 6\\
	73 12\\
	73 15\\
	74 6\\
	74 7\\
	74 15\\
	74 16\\
	75 7\\
	75 8\\
	75 16\\
	75 17\\
	76 2\\
	76 8\\
	76 11\\
	76 17\\
	77 8\\
	77 9\\
	77 17\\
	77 18\\
	78 4\\
	78 9\\
	78 13\\
	78 18\\
	1 43\\
	4 43\\
	10 43\\
	13 43\\
	19 43\\
	2 44\\
	8 44\\
	11 44\\
	17 44\\
	20 44\\
	3 45\\
	6 45\\
	12 45\\
	15 45\\
	21 45\\
	1 46\\
	4 46\\
	5 46\\
	9 46\\
	10 46\\
	13 46\\
	14 46\\
	18 46\\
	4 47\\
	5 47\\
	6 47\\
	13 47\\
	14 47\\
	15 47\\
	3 48\\
	5 48\\
	6 48\\
	7 48\\
	12 48\\
	14 48\\
	15 48\\
	16 48\\
	6 49\\
	7 49\\
	8 49\\
	15 49\\
	16 49\\
	17 49\\
	2 50\\
	7 50\\
	8 50\\
	9 50\\
	11 50\\
	16 50\\
	17 50\\
	18 50\\
	4 51\\
	8 51\\
	9 51\\
	13 51\\
	17 51\\
	18 51\\
	1 52\\
	4 52\\
	10 52\\
	13 52\\
	22 52\\
	2 53\\
	8 53\\
	11 53\\
	17 53\\
	23 53\\
	3 54\\
	6 54\\
	12 54\\
	15 54\\
	24 54\\
	1 55\\
	4 55\\
	5 55\\
	9 55\\
	10 55\\
	13 55\\
	14 55\\
	18 55\\
	4 56\\
	5 56\\
	6 56\\
	13 56\\
	14 56\\
	15 56\\
	3 57\\
	5 57\\
	6 57\\
	7 57\\
	12 57\\
	14 57\\
	15 57\\
	16 57\\
	6 58\\
	7 58\\
	8 58\\
	15 58\\
	16 58\\
	17 58\\
	2 59\\
	7 59\\
	8 59\\
	9 59\\
	11 59\\
	16 59\\
	17 59\\
	18 59\\
	4 60\\
	8 60\\
	9 60\\
	13 60\\
	17 60\\
	18 60\\
	1 61\\
	4 61\\
	10 61\\
	13 61\\
	4 62\\
	5 62\\
	13 62\\
	14 62\\
	5 63\\
	6 63\\
	14 63\\
	15 63\\
	3 64\\
	6 64\\
	12 64\\
	15 64\\
	6 65\\
	7 65\\
	15 65\\
	16 65\\
	7 66\\
	8 66\\
	16 66\\
	17 66\\
	2 67\\
	8 67\\
	11 67\\
	17 67\\
	8 68\\
	9 68\\
	17 68\\
	18 68\\
	4 69\\
	9 69\\
	13 69\\
	18 69\\
	1 70\\
	4 70\\
	10 70\\
	13 70\\
	4 71\\
	5 71\\
	13 71\\
	14 71\\
	5 72\\
	6 72\\
	14 72\\
	15 72\\
	3 73\\
	6 73\\
	12 73\\
	15 73\\
	6 74\\
	7 74\\
	15 74\\
	16 74\\
	7 75\\
	8 75\\
	16 75\\
	17 75\\
	2 76\\
	8 76\\
	11 76\\
	17 76\\
	8 77\\
	9 77\\
	17 77\\
	18 77\\
	4 78\\
	9 78\\
	13 78\\
	18 78\\
	    43    43\\
	44    44\\
	45    45\\
	46    46\\
	47    47\\
	48    48\\
	49    49\\
	50    50\\
	51    51\\
	52    52\\
	53    53\\
	54    54\\
	55    55\\
	56    56\\
	57    57\\
	58    58\\
	59    59\\
	60    60\\
	61    61\\
	62    62\\
	63    63\\
	64    64\\
	65    65\\
	66    66\\
	67    67\\
	68    68\\
	69    69\\
	70    70\\
	71    71\\
	72    72\\
	73    73\\
	74    74\\
	75    75\\
	76    76\\
	77    77\\
	78    78\\
};

\addplot [color=black, draw=none, mark size=1pt, mark=*, mark options={solid, mycolor3}]
table[row sep=crcr]{%
	    25 25\\
	26 26\\
	27 27\\
	28 28\\
	29 29\\
	30 30\\
	31 31\\
	32 32\\
	33 33\\
	34 34\\
	35 35\\
	36 36\\
	37 37\\
	38 38\\
	39 39\\
	40 40\\
	41 41\\
	42 42\\
	61 25\\
	62 26\\
	63 27\\
	64 28\\
	65 29\\
	66 30\\
	67 31\\
	68 32\\
	69 33\\
	70 34\\
	71 35\\
	72 36\\
	73 37\\
	74 38\\
	75 39\\
	76 40\\
	77 41\\
	78 42\\
	25 25\\
	26 26\\
	27 27\\
	28 28\\
	29 29\\
	30 30\\
	31 31\\
	32 32\\
	33 33\\
	34 34\\
	35 35\\
	36 36\\
	37 37\\
	38 38\\
	39 39\\
	40 40\\
	41 41\\
	42 42\\
	25 61\\
	26 62\\
	27 63\\
	28 64\\
	29 65\\
	30 66\\
	31 67\\
	32 68\\
	33 69\\
	34 70\\
	35 71\\
	36 72\\
	37 73\\
	38 74\\
	39 75\\
	40 76\\
	41 77\\
	42 78\\
};

\end{axis}
\end{tikzpicture}%
  ~
%
%
\begin{tikzpicture}[scale=.55]

\begin{axis}[%
width=0.5\textwidth,
height=0.5\textwidth,
at={(1.648in,0.642in)},
scale only axis,
xmin=0,
xmax=79,
xlabel style={font=\color{white!15!black}},
xlabel={nz = 576},
y dir=reverse,
ymin=0,
ymax=79,
axis background/.style={fill=white},
legend style={legend cell align=left, align=left, draw=white!15!black},
every x tick label/.append style={font=\normalsize},
every y tick label/.append style={font=\normalsize},
label style={font=\large}
]

\addplot [color=black, draw=none, mark size=1pt, mark=*, mark options={solid, mycolor0}]
table[row sep=crcr]{%
	1 1\\
	1 4\\
	1 10\\
	1 13\\
	2 2\\
	2 8\\
	2 11\\
	2 17\\
	3 3\\
	3 6\\
	3 12\\
	3 15\\
	4 1\\
	4 4\\
	4 5\\
	4 9\\
	4 10\\
	4 13\\
	4 14\\
	4 18\\
	5 4\\
	5 5\\
	5 6\\
	5 13\\
	5 14\\
	5 15\\
	6 3\\
	6 5\\
	6 6\\
	6 7\\
	6 12\\
	6 14\\
	6 15\\
	6 16\\
	7 6\\
	7 7\\
	7 8\\
	7 15\\
	7 16\\
	7 17\\
	8 2\\
	8 7\\
	8 8\\
	8 9\\
	8 11\\
	8 16\\
	8 17\\
	8 18\\
	9 4\\
	9 8\\
	9 9\\
	9 13\\
	9 17\\
	9 18\\
	10 1\\
	10 4\\
	10 10\\
	10 13\\
	11 2\\
	11 8\\
	11 11\\
	11 17\\
	12 3\\
	12 6\\
	12 12\\
	12 15\\
	13 1\\
	13 4\\
	13 5\\
	13 9\\
	13 10\\
	13 13\\
	13 14\\
	13 18\\
	14 4\\
	14 5\\
	14 6\\
	14 13\\
	14 14\\
	14 15\\
	15 3\\
	15 5\\
	15 6\\
	15 7\\
	15 12\\
	15 14\\
	15 15\\
	15 16\\
	16 6\\
	16 7\\
	16 8\\
	16 15\\
	16 16\\
	16 17\\
	17 2\\
	17 7\\
	17 8\\
	17 9\\
	17 11\\
	17 16\\
	17 17\\
	17 18\\
	18 4\\
	18 8\\
	18 9\\
	18 13\\
	18 17\\
	18 18\\
	19 19\\
	20 20\\
	21 21\\
	22 22\\
	23 23\\
	24 24\\
};

\addplot [color=black, draw=none, mark size=1pt, mark=*, mark options={solid, mycolor1}]
table[row sep=crcr]{%
	    24 1\\
	24 4\\
	24 10\\
	24 13\\
	24 19\\
	25 2\\
	25 8\\
	25 11\\
	25 17\\
	25 20\\
	26 3\\
	26 6\\
	26 12\\
	26 15\\
	26 21\\
	27 1\\
	27 4\\
	27 5\\
	27 9\\
	27 10\\
	27 13\\
	27 14\\
	27 18\\
	28 4\\
	28 5\\
	28 6\\
	28 13\\
	28 14\\
	28 15\\
	29 3\\
	29 5\\
	29 6\\
	29 7\\
	29 12\\
	29 14\\
	29 15\\
	29 16\\
	30 6\\
	30 7\\
	30 8\\
	30 15\\
	30 16\\
	30 17\\
	31 2\\
	31 7\\
	31 8\\
	31 9\\
	31 11\\
	31 16\\
	31 17\\
	31 18\\
	32 4\\
	32 8\\
	32 9\\
	32 13\\
	32 17\\
	32 18\\
	33 1\\
	33 4\\
	33 10\\
	33 13\\
	33 22\\
	34 2\\
	34 8\\
	34 11\\
	34 17\\
	34 23\\
	35 3\\
	35 6\\
	35 12\\
	35 15\\
	35 24\\
	36 1\\
	36 4\\
	36 5\\
	36 9\\
	36 10\\
	36 13\\
	36 14\\
	36 18\\
	37 4\\
	37 5\\
	37 6\\
	37 13\\
	37 14\\
	37 15\\
	38 3\\
	38 5\\
	38 6\\
	38 7\\
	38 12\\
	38 14\\
	38 15\\
	38 16\\
	39 6\\
	39 7\\
	39 8\\
	39 15\\
	39 16\\
	39 17\\
	40 2\\
	40 7\\
	40 8\\
	40 9\\
	40 11\\
	40 16\\
	40 17\\
	40 18\\
	41 4\\
	41 8\\
	41 9\\
	41 13\\
	41 17\\
	41 18\\
	42 1\\
	42 4\\
	42 10\\
	42 13\\
	43 4\\
	43 5\\
	43 13\\
	43 14\\
	44 5\\
	44 6\\
	44 14\\
	44 15\\
	45 3\\
	45 6\\
	45 12\\
	45 15\\
	46 6\\
	46 7\\
	46 15\\
	46 16\\
	47 7\\
	47 8\\
	47 16\\
	47 17\\
	48 2\\
	48 8\\
	48 11\\
	48 17\\
	49 8\\
	49 9\\
	49 17\\
	49 18\\
	50 4\\
	50 9\\
	50 13\\
	50 18\\
	51 1\\
	51 4\\
	51 10\\
	51 13\\
	52 4\\
	52 5\\
	52 13\\
	52 14\\
	53 5\\
	53 6\\
	53 14\\
	53 15\\
	54 3\\
	54 6\\
	54 12\\
	54 15\\
	55 6\\
	55 7\\
	55 15\\
	55 16\\
	56 7\\
	56 8\\
	56 16\\
	56 17\\
	57 2\\
	57 8\\
	57 11\\
	57 17\\
	58 8\\
	58 9\\
	58 17\\
	58 18\\
	59 4\\
	59 9\\
	59 13\\
	59 18\\
	1 24\\
	4 24\\
	10 24\\
	13 24\\
	19 24\\
	2 25\\
	8 25\\
	11 25\\
	17 25\\
	20 25\\
	3 26\\
	6 26\\
	12 26\\
	15 26\\
	21 26\\
	1 27\\
	4 27\\
	5 27\\
	9 27\\
	10 27\\
	13 27\\
	14 27\\
	18 27\\
	4 28\\
	5 28\\
	6 28\\
	13 28\\
	14 28\\
	15 28\\
	3 29\\
	5 29\\
	6 29\\
	7 29\\
	12 29\\
	14 29\\
	15 29\\
	16 29\\
	6 30\\
	7 30\\
	8 30\\
	15 30\\
	16 30\\
	17 30\\
	2 31\\
	7 31\\
	8 31\\
	9 31\\
	11 31\\
	16 31\\
	17 31\\
	18 31\\
	4 32\\
	8 32\\
	9 32\\
	13 32\\
	17 32\\
	18 32\\
	1 33\\
	4 33\\
	10 33\\
	13 33\\
	22 33\\
	2 34\\
	8 34\\
	11 34\\
	17 34\\
	23 34\\
	3 35\\
	6 35\\
	12 35\\
	15 35\\
	24 35\\
	1 36\\
	4 36\\
	5 36\\
	9 36\\
	10 36\\
	13 36\\
	14 36\\
	18 36\\
	4 37\\
	5 37\\
	6 37\\
	13 37\\
	14 37\\
	15 37\\
	3 38\\
	5 38\\
	6 38\\
	7 38\\
	12 38\\
	14 38\\
	15 38\\
	16 38\\
	6 39\\
	7 39\\
	8 39\\
	15 39\\
	16 39\\
	17 39\\
	2 40\\
	7 40\\
	8 40\\
	9 40\\
	11 40\\
	16 40\\
	17 40\\
	18 40\\
	4 41\\
	8 41\\
	9 41\\
	13 41\\
	17 41\\
	18 41\\
	1 42\\
	4 42\\
	10 42\\
	13 42\\
	4 43\\
	5 43\\
	13 43\\
	14 43\\
	5 44\\
	6 44\\
	14 44\\
	15 44\\
	3 45\\
	6 45\\
	12 45\\
	15 45\\
	6 46\\
	7 46\\
	15 46\\
	16 46\\
	7 47\\
	8 47\\
	16 47\\
	17 47\\
	2 48\\
	8 48\\
	11 48\\
	17 48\\
	8 49\\
	9 49\\
	17 49\\
	18 49\\
	4 50\\
	9 50\\
	13 50\\
	18 50\\
	1 51\\
	4 51\\
	10 51\\
	13 51\\
	4 52\\
	5 52\\
	13 52\\
	14 52\\
	5 53\\
	6 53\\
	14 53\\
	15 53\\
	3 54\\
	6 54\\
	12 54\\
	15 54\\
	6 55\\
	7 55\\
	15 55\\
	16 55\\
	7 56\\
	8 56\\
	16 56\\
	17 56\\
	2 57\\
	8 57\\
	11 57\\
	17 57\\
	8 58\\
	9 58\\
	17 58\\
	18 58\\
	4 59\\
	9 59\\
	13 59\\
	18 59\\
	    24    24\\
	25    25\\
	26    26\\
	27    27\\
	28    28\\
	29    29\\
	30    30\\
	31    31\\
	32    32\\
	33    33\\
	34    34\\
	35    35\\
	36    36\\
	37    37\\
	38    38\\
	39    39\\
	40    40\\
	41    41\\
	42    42\\
	43    43\\
	44    44\\
	45    45\\
	46    46\\
	47    47\\
	48    48\\
	49    49\\
	50    50\\
	51    51\\
	52    52\\
	53    53\\
	54    54\\
	55    55\\
	56    56\\
	57    57\\
	58    58\\
	59    59\\
	60    60\\
};

\addplot [color=black, draw=none, mark size=1pt, mark=*, mark options={solid, mycolor3}]
table[row sep=crcr]{%
	    61 43\\
	62 44\\
	63 45\\
	64 46\\
	65 47\\
	66 48\\
	67 49\\
	68 50\\
	69 51\\
	70 52\\
	71 53\\
	72 54\\
	73 55\\
	74 56\\
	75 57\\
	76 58\\
	77 59\\
	78 60\\
	    43 61\\
	44 62\\
	45 63\\
	46 64\\
	47 65\\
	48 66\\
	49 67\\
	50 68\\
	51 69\\
	52 70\\
	53 71\\
	54 72\\
	55 73\\
	56 74\\
	57 75\\
	58 76\\
	59 77\\
	60 78\\
	    61    61\\
	62    62\\
	63    63\\
	64    64\\
	65    65\\
	66    66\\
	67    67\\
	68    68\\
	69    69\\
	70    70\\
	71    71\\
	72    72\\
	73    73\\
	74    74\\
	75    75\\
	76    76\\
	77    77\\
	78    78\\
};

\end{axis}
\end{tikzpicture}%
  ~ 
%
%
\begin{tikzpicture}[scale=.55]

\begin{axis}[%
width=0.5\textwidth,
height=0.5\textwidth,
at={(1.648in,0.642in)},
scale only axis,
xmin=0,
xmax=61,
xlabel style={font=\color{white!15!black}},
xlabel={nz = 522},
y dir=reverse,
ymin=0,
ymax=61,
axis background/.style={fill=white},
legend style={legend cell align=left, align=left, draw=white!15!black},
every x tick label/.append style={font=\normalsize},
every y tick label/.append style={font=\normalsize},
label style={font=\large}
]

\addplot [color=black, draw=none, mark size=1pt, mark=*, mark options={solid, mycolor0}]
table[row sep=crcr]{%
	1 1\\
	1 4\\
	1 10\\
	1 13\\
	2 2\\
	2 8\\
	2 11\\
	2 17\\
	3 3\\
	3 6\\
	3 12\\
	3 15\\
	4 1\\
	4 4\\
	4 5\\
	4 9\\
	4 10\\
	4 13\\
	4 14\\
	4 18\\
	5 4\\
	5 5\\
	5 6\\
	5 13\\
	5 14\\
	5 15\\
	6 3\\
	6 5\\
	6 6\\
	6 7\\
	6 12\\
	6 14\\
	6 15\\
	6 16\\
	7 6\\
	7 7\\
	7 8\\
	7 15\\
	7 16\\
	7 17\\
	8 2\\
	8 7\\
	8 8\\
	8 9\\
	8 11\\
	8 16\\
	8 17\\
	8 18\\
	9 4\\
	9 8\\
	9 9\\
	9 13\\
	9 17\\
	9 18\\
	10 1\\
	10 4\\
	10 10\\
	10 13\\
	11 2\\
	11 8\\
	11 11\\
	11 17\\
	12 3\\
	12 6\\
	12 12\\
	12 15\\
	13 1\\
	13 4\\
	13 5\\
	13 9\\
	13 10\\
	13 13\\
	13 14\\
	13 18\\
	14 4\\
	14 5\\
	14 6\\
	14 13\\
	14 14\\
	14 15\\
	15 3\\
	15 5\\
	15 6\\
	15 7\\
	15 12\\
	15 14\\
	15 15\\
	15 16\\
	16 6\\
	16 7\\
	16 8\\
	16 15\\
	16 16\\
	16 17\\
	17 2\\
	17 7\\
	17 8\\
	17 9\\
	17 11\\
	17 16\\
	17 17\\
	17 18\\
	18 4\\
	18 8\\
	18 9\\
	18 13\\
	18 17\\
	18 18\\
	19 19\\
	20 20\\
	21 21\\
	22 22\\
	23 23\\
	24 24\\
};

\addplot [color=black, draw=none, mark size=1pt, mark=*, mark options={solid, mycolor1}]
table[row sep=crcr]{%
	24 1\\
	24 4\\
	24 10\\
	24 13\\
	24 19\\
	25 2\\
	25 8\\
	25 11\\
	25 17\\
	25 20\\
	26 3\\
	26 6\\
	26 12\\
	26 15\\
	26 21\\
	27 1\\
	27 4\\
	27 5\\
	27 9\\
	27 10\\
	27 13\\
	27 14\\
	27 18\\
	28 4\\
	28 5\\
	28 6\\
	28 13\\
	28 14\\
	28 15\\
	29 3\\
	29 5\\
	29 6\\
	29 7\\
	29 12\\
	29 14\\
	29 15\\
	29 16\\
	30 6\\
	30 7\\
	30 8\\
	30 15\\
	30 16\\
	30 17\\
	31 2\\
	31 7\\
	31 8\\
	31 9\\
	31 11\\
	31 16\\
	31 17\\
	31 18\\
	32 4\\
	32 8\\
	32 9\\
	32 13\\
	32 17\\
	32 18\\
	33 1\\
	33 4\\
	33 10\\
	33 13\\
	33 22\\
	34 2\\
	34 8\\
	34 11\\
	34 17\\
	34 23\\
	35 3\\
	35 6\\
	35 12\\
	35 15\\
	35 24\\
	36 1\\
	36 4\\
	36 5\\
	36 9\\
	36 10\\
	36 13\\
	36 14\\
	36 18\\
	37 4\\
	37 5\\
	37 6\\
	37 13\\
	37 14\\
	37 15\\
	38 3\\
	38 5\\
	38 6\\
	38 7\\
	38 12\\
	38 14\\
	38 15\\
	38 16\\
	39 6\\
	39 7\\
	39 8\\
	39 15\\
	39 16\\
	39 17\\
	40 2\\
	40 7\\
	40 8\\
	40 9\\
	40 11\\
	40 16\\
	40 17\\
	40 18\\
	41 4\\
	41 8\\
	41 9\\
	41 13\\
	41 17\\
	41 18\\
	42 1\\
	42 4\\
	42 10\\
	42 13\\
	43 4\\
	43 5\\
	43 13\\
	43 14\\
	44 5\\
	44 6\\
	44 14\\
	44 15\\
	45 3\\
	45 6\\
	45 12\\
	45 15\\
	46 6\\
	46 7\\
	46 15\\
	46 16\\
	47 7\\
	47 8\\
	47 16\\
	47 17\\
	48 2\\
	48 8\\
	48 11\\
	48 17\\
	49 8\\
	49 9\\
	49 17\\
	49 18\\
	50 4\\
	50 9\\
	50 13\\
	50 18\\
	51 1\\
	51 4\\
	51 10\\
	51 13\\
	52 4\\
	52 5\\
	52 13\\
	52 14\\
	53 5\\
	53 6\\
	53 14\\
	53 15\\
	54 3\\
	54 6\\
	54 12\\
	54 15\\
	55 6\\
	55 7\\
	55 15\\
	55 16\\
	56 7\\
	56 8\\
	56 16\\
	56 17\\
	57 2\\
	57 8\\
	57 11\\
	57 17\\
	58 8\\
	58 9\\
	58 17\\
	58 18\\
	59 4\\
	59 9\\
	59 13\\
	59 18\\
	1 24\\
	4 24\\
	10 24\\
	13 24\\
	19 24\\
	2 25\\
	8 25\\
	11 25\\
	17 25\\
	20 25\\
	3 26\\
	6 26\\
	12 26\\
	15 26\\
	21 26\\
	1 27\\
	4 27\\
	5 27\\
	9 27\\
	10 27\\
	13 27\\
	14 27\\
	18 27\\
	4 28\\
	5 28\\
	6 28\\
	13 28\\
	14 28\\
	15 28\\
	3 29\\
	5 29\\
	6 29\\
	7 29\\
	12 29\\
	14 29\\
	15 29\\
	16 29\\
	6 30\\
	7 30\\
	8 30\\
	15 30\\
	16 30\\
	17 30\\
	2 31\\
	7 31\\
	8 31\\
	9 31\\
	11 31\\
	16 31\\
	17 31\\
	18 31\\
	4 32\\
	8 32\\
	9 32\\
	13 32\\
	17 32\\
	18 32\\
	1 33\\
	4 33\\
	10 33\\
	13 33\\
	22 33\\
	2 34\\
	8 34\\
	11 34\\
	17 34\\
	23 34\\
	3 35\\
	6 35\\
	12 35\\
	15 35\\
	24 35\\
	1 36\\
	4 36\\
	5 36\\
	9 36\\
	10 36\\
	13 36\\
	14 36\\
	18 36\\
	4 37\\
	5 37\\
	6 37\\
	13 37\\
	14 37\\
	15 37\\
	3 38\\
	5 38\\
	6 38\\
	7 38\\
	12 38\\
	14 38\\
	15 38\\
	16 38\\
	6 39\\
	7 39\\
	8 39\\
	15 39\\
	16 39\\
	17 39\\
	2 40\\
	7 40\\
	8 40\\
	9 40\\
	11 40\\
	16 40\\
	17 40\\
	18 40\\
	4 41\\
	8 41\\
	9 41\\
	13 41\\
	17 41\\
	18 41\\
	1 42\\
	4 42\\
	10 42\\
	13 42\\
	4 43\\
	5 43\\
	13 43\\
	14 43\\
	5 44\\
	6 44\\
	14 44\\
	15 44\\
	3 45\\
	6 45\\
	12 45\\
	15 45\\
	6 46\\
	7 46\\
	15 46\\
	16 46\\
	7 47\\
	8 47\\
	16 47\\
	17 47\\
	2 48\\
	8 48\\
	11 48\\
	17 48\\
	8 49\\
	9 49\\
	17 49\\
	18 49\\
	4 50\\
	9 50\\
	13 50\\
	18 50\\
	1 51\\
	4 51\\
	10 51\\
	13 51\\
	4 52\\
	5 52\\
	13 52\\
	14 52\\
	5 53\\
	6 53\\
	14 53\\
	15 53\\
	3 54\\
	6 54\\
	12 54\\
	15 54\\
	6 55\\
	7 55\\
	15 55\\
	16 55\\
	7 56\\
	8 56\\
	16 56\\
	17 56\\
	2 57\\
	8 57\\
	11 57\\
	17 57\\
	8 58\\
	9 58\\
	17 58\\
	18 58\\
	4 59\\
	9 59\\
	13 59\\
	18 59\\
	24    24\\
	25    25\\
	26    26\\
	27    27\\
	28    28\\
	29    29\\
	30    30\\
	31    31\\
	32    32\\
	33    33\\
	34    34\\
	35    35\\
	36    36\\
	37    37\\
	38    38\\
	39    39\\
	40    40\\
	41    41\\
	42    42\\
	43    43\\
	44    44\\
	45    45\\
	46    46\\
	47    47\\
	48    48\\
	49    49\\
	50    50\\
	51    51\\
	52    52\\
	53    53\\
	54    54\\
	55    55\\
	56    56\\
	57    57\\
	58    58\\
	59    59\\
	60    60\\
};
\end{axis}
\end{tikzpicture}%
 \caption{Structure of the KKT system \eqref{eq:SCOPFhessianSymmetric}, reordered according to \eqref{eq:SCOPFslackAtEnd}, and the structure of the reduced KKT with the slacks removed \eqref{eq:SCOPFhessianReduced}.}
\label{fig:KKTsystemOPF}
\end{center}
\end{figure}

For a robust algorithm it is crucial to obtain highly accurate search directions. Most of the burden is shifted to the sparse linear solver, where techniques such as fill-in minimization reordering, symmetric scaling vectors, matching, and pivoting can provide substantial improvement to the solution accuracy. Additional improvement can be achieved by performing iterative refinement using the unsymmetrical version KKT linear system of form \eqref{eq:SCOPFhessian}. 
It is possible to further reduce the KKT system by eliminating the slack variables $\bfs$. The system \eqref{eq:SCOPFhessianSymmetric} can be permuted to the structure with the diagonal block $L_s$ in the lower right corner,
\begin{equation}
\label{eq:SCOPFslackAtEnd}
\begin{bmatrix}
\tilde{\Hess} & \Je^\transpose & \Ji^\transpose & \0 \\
\Je & \0 & \0 & \0 & \\
\Ji & \0 & \0 & -I \\
\0 & \0 & -I & L_s
\end{bmatrix}^k
\begin{bmatrix}
\Delta \bfx\\
\Delta \lamE\\
\Delta \lamI\\
\Delta \bfs
\end{bmatrix}^k
=
-\begin{bmatrix}
\boldsymbol{l}_a + {X}^{-1} \boldsymbol{l}_e\\
\boldsymbol{l}_c \\
\boldsymbol{l}_d\\
\boldsymbol{l}_b + {S}^{-1} \boldsymbol{l}_f
\end{bmatrix}^k.
\end{equation}
Since the block $L_s$ is a diagonal matrix, the reordered system \eqref{eq:SCOPFslackAtEnd} can be trivially reduced by computing the Schur complement with respect to the $3 \times 3$ block in the upper left corner, as illustrated in Figure \ref{fig:KKTsystemOPF},
\begin{equation}
\label{eq:SCOPFhessianSchurComplement}
\begin{bmatrix}
\tilde{\Hess} & \Je^\transpose & \Ji^\transpose \\
\Je & \0 & \0 &  \\
\Ji & \0 & \0
\end{bmatrix}^k -
\begin{bmatrix}
 \0 & \0 & -I
\end{bmatrix}^\transpose {(L_s^k)}^{-1}
\begin{bmatrix}
 \0 & \0 & -I
\end{bmatrix}.
\end{equation}
The additional elimination, compared to \cite{olafAugmented, pardiso}, further reduces the memory requirements and computation time due to the smaller amount of factorization fill-in. Such an elimination, however, can be performed only for the nonzero elements of $L_s^k$ sufficiently away from zero in order to avoid the ill-conditioning of the reduced system.
The reduced linear system that needs to be solved now has the structure
\begin{equation}
\label{eq:SCOPFhessianReduced}
\begin{bmatrix}
\tilde{\Hess} & \Je^\transpose & \Ji^\transpose \\
\Je & \0 & \0 &  \\
\Ji & \0 & - L_s^{-1}
\end{bmatrix}^k
\begin{bmatrix}
\Delta \bfx\\
\Delta \lamE\\
\Delta \lamI\\
\end{bmatrix}^k
=
-\begin{bmatrix}
\boldsymbol{l}_a + {X}^{-1} \boldsymbol{l}_e\\
\boldsymbol{l}_c \\
\boldsymbol{l}_d + L_s^{-1}(\boldsymbol{l}_b + S^{-1}\boldsymbol{l}_f) \\
\end{bmatrix}^k
\end{equation}
and the eliminated slack variables can be recovered by solving
\begin{equation}
\label{eq:SCOPFslackReduction}
L_s^k \Delta \bfs^k = - \boldsymbol{l}_b^k - S_k^{-1} \boldsymbol{l}_f^k + \Delta \lamI^k.
\end{equation}



\subsection{Backtracking Line-Search Filter Method  \label{sec:FilterLS}}
After the successful computation of the search direction from \eqref{eq:SCOPFhessianSymmetric} and \eqref{eq:SCOPFhessianSymmetricRecovery} the step sizes $\alpha_k, \alpha_k^z \in (0,1]$ need to be determined in order to obtain the next iterate:
\begingroup
\allowdisplaybreaks
\begin{align}
\bfx^{k+1} &= \bfx^{k} + \alpha_k \Delta \bfx^k,\\
\bfs^{k+1} &= \bfs^{k} + \alpha_k \Delta \bfs^k,\\
\lamE^{k+1} &= \lamE^{k} + \alpha_k \Delta \lamE^k,\\
\lamI^{k+1} &= \lamI^{k} + \alpha_k \Delta \lamI^k,\\
\bfz^{k+1} &= \bfz^{k} + \alpha_k^z \Delta \bfz^k, \\
\bfy^{k+1} &= \bfy^{k} + \alpha_k^z \Delta \bfy^k.
\end{align}
\endgroup
Different step sizes for the primal and dual variables is commonly employed to prevent unnecessarily small steps in either variables and delay the convergence to the optimal. A first candidate step length is chosen such that the strict positivity of $\bfx, \bfs$, and $\bfz$ is preserved, since it needs to hold both in the solution of the barrier problem \eqref{eq:Barrier} and also in every iteration, which is necessary in order to evaluate the barrier function. This is accomplished by the fraction-to-the-boundary rule, which identifies the maximum step size $\alpha_k, \alpha^z_k \in (0,1]$, such that
\begin{align}
    \alpha^{\text{max}}_k &= \max \left(\alpha \in (0,1]: \bfx^{k} + \alpha \Delta \bfx^k \geq (1-\tau)\bfx^k \right),\\
    \alpha_k^z &= \max \left(\alpha \in (0,1]: \bfz^{k} + \alpha \Delta \bfz^k \geq (1-\tau)\bfz^k \right),
\end{align}
where $\tau \in (0,1)$ is a function of the current barrier parameter $\mu_j$.  The step size for the dual variables $\alpha^z_k$ is used directly, but in order to ensure global convergence the step size $\alpha_k \in (0,\alpha_k^{\text{max}})$ for the remaining variables is determined by a backtracking line-search procedure, exploring a decreasing sequence of trial step sizes $\alpha_k^i = 2^{-i} \alpha_k^{\text{max}}$ for $i = 0,1,2,\ldots$.

The variant of the backtracking line-search filter method \cite{Fletcher2002} used in \IPOPT{} is based on the idea of a biobjective optimization problem with the two goals (i) minimizing the objective function 
\begin{equation}
\varphi_{\mu_j}(\bfx,\bfs) := f(\bfx) -  \mu_j \sum_{i=1}^{n} \log(x_i) - \mu_j \sum_{i=1}^{\ni} \log(s_i),
\end{equation}
and (ii) minimizing the constraint violation
\begin{equation}
\theta(\bfx,\bfs) := \Vert\  (\ce(\bfx),\ \ci(\bfx)-\bfs) \  \Vert_1.    
\end{equation}
A trial point
$\bfx^{k}(\alpha_k^i) := \bfx^{k} + \alpha_k^i \Delta \bfx^k$ and 
$\bfs^{k}(\alpha_k^i) := \bfs^{k} + \alpha_k^i \Delta \bfs^k$
during the backtracking line search is considered to be acceptable, if it leads to sufficient progress toward either goal compared to the current iterate. The emphasis is put on the latter goal, until the constraint violations satisfy a certain threshold. Afterwards, the former goal is emphasized and reduction in the barrier function is required, accepting only iterates satisfying the Armijo condition.

\begin{definition}
The filter $\mathcal{F}$ is a set of ordered pairs containing a constraint violation value $\theta$ and the objective function value $\varphi$,
such that
\begin{equation}
    \mathcal{F} \subseteq \{(\theta, \varphi) \in \mathbb{R}^2 : \theta >0 \}.
\end{equation}
\end{definition}
The algorithm also maintains a filter $\mathcal{F}_j$ for each BSP $j$ for which the $\mu_j$ is fixed. The filter $\mathcal{F}_j$ contains those combinations that are prohibited for a successful trial point in all iterations within the $j$th BSP. The filter is initialized so that the algorithm will never allow trial points to be accepted that have a constraint violation larger than $\theta^{\text{max}}$. During the line search, a trial point $\bfx^k(\alpha_k^i),\ \bfs^{k}(\alpha_k^i)$ is rejected if $(\theta(\bfx_k(\alpha_k^i), \bfs^{k}(\alpha_k^i)),\ \varphi_{\mu_j}(\bfx^k(\alpha_k^i), \bfs^{k}(\alpha_k^i))) \in \mathcal{F}_j$. After every iteration, in which the accepted trial step size does not satisfy the two objectives of the backtracking linesearch, the filter is augmented. This ensures that the iterates cannot return to the neighborhood of the unsatisfactory iterates. Overall, this procedure ensures that the algorithm cannot cycle, for example, between two points that alternate between decrease of the constraint violation and the barrier objective function.

In cases when it is not possible to identify a satisfactory trial step size, the algorithm reverts to a feasibility restoration phase. Here, the algorithm tries to find a new iterate which is acceptable to the current filter, by reducing the constraint violation with some iterative method. Note that the restoration phase algorithm might not be able to produce a new iterate for the filter line-search method, for example, when the problem is infeasible.
\subsection{Inertia Correction and Curvature Detection \label{sec:inertia}}
\begin{definition}
The inertia of a square matrix is defined as the ordered triplet $(n+, n-, n_0) \in \{\mathbb{N}\cup 0\}^3$, where the terms denote the number of positive, negative, and zero eigenvalues, respectively.
\end{definition}

In order to guarantee descent properties for the line-search procedure, it is necessary to ensure that the Hessian matrix projected on the null space of the constraint Jacobian is positive definite (see Theorem \ref{theorem:SOC}). Also, if the constraint Jacobian does not have full rank, the iteration matrix in \eqref{eq:SCOPFhessianSymmetric} is singular, and the solution might not exist. These conditions are satisfied if the iteration matrix has the inertia $(N_x+\ni,\ne+\ni,0)$. The sizes correspond to the size of the Hessian block (with respect to both primal variables $\bfx$ and the slack variables $\bfs$) and the Jacobians of the equality and inequality constraints. If the inertia is not correct, the iteration matrix needs to be modified. In \IPOPT{} implementation, the diagonal perturbations $\delta_w, \delta_c \geq 0$ are added to the Hessian \eqref{eq:SCOPFhessianSymmetric}, such that
\begin{equation}
    \begin{bmatrix}
    \boldsymbol{\tilde{H}} + \delta_w I & \0 & \Je^\transpose & \Ji^\transpose \\
    \0 & {L}_{s} + \delta_w I & \0 & -I \\
    \Je & \0 & -\delta_c I & \0 & \\
    \Ji & -I & \0 & -\delta_c I & 
    \end{bmatrix}.
\end{equation}
The system is refactorized with different trial values of $\delta_w, \delta_c$ until the inertia is correct. The inertia of the iteration matrix is readily available from several sparse indefinite linear solvers, such as \PARDISO{} \cite{SCHENK2001}. In case the correct inertia cannot be achieved, the current search direction computation is aborted and the algorithm uses a different objective function that does try to solely minimize the feasibility violation (e.g., minimizing the constraints violation), ignoring the original objective function, in the hope that the matrix has better properties close to the feasible points.

The inertia detection strategy focuses on the properties of the augmented iteration matrix \eqref{eq:SCOPFhessianSymmetric} alone and can discard search directions that are of descent but for which the inertia of the augmented matrix is not correct. Furthermore, the inertia detection strategy might require multiple factorizations of the iteration matrix and, because the factorization is the most expensive step in the algorithm, computational performance can be greatly affected. Furthermore, the inertia estimates might vary, depending on which linear solver is used or not be available at all. To bypass the need for the inertia information, several authors suggest using the curvature test, e.g., \cite{CosminCurvature} \cite{Chiang2016}:
\begin{align}
\label{eq:curvature}
& \bfd_k^\transpose \boldsymbol{W}_k(\delta) \, \bfd_k \geq \kappa \bfd_k^\transpose \bfd_k, \quad \kappa > 0, \delta \geq 0,\\
\nonumber
& \boldsymbol{W}_k(\delta) = 
 \begin{bmatrix}
 \boldsymbol{\tilde{H}} & \0 \\
 \0 & {L}_{s}
 \end{bmatrix}^k + \delta I, \quad \bfd_k = \pmat{\Delta \bfx_k, & \Delta \bfs_k}.
\end{align}
If the test is satisfied, the search direction is accepted; if it is not satisfied, the regularization parameter $\delta$ is increased and a new search direction is computed using the new regularized matrix. 

\begin{remark}
While the curvature detection strategy usually requires more IP iterations until convergence compared with the inertia detection, it may require fewer extra factorizations. Overall, the solution time is less than that of the inertia detection because significantly fewer regularizations are needed. 
\end{remark}

\subsection{Barrier Parameter Update Strategy  \label{sec:Barrier}}
The strategy of the barrier parameter update is an important factor influencing the convergence properties, especially for difficult nonconvex problems.
When solving nonlinear nonconvex programming problems, it is of great importance to prevent the iteration from failing. 
Different barrier parameter update strategies are discussed here, including the monotone Fiacco--McCormick strategy \cite{muMONOTONE} and an adaptive strategy based on minimization of a quality function \cite{muADAPTIVEqualityfcn}.


\subsubsection{Monotone and Adaptive Strategies}
Using the default monotone Fiacco--McCormick strategy, an approximate solution to the barrier problem \eqref{eq:Barrier} for a fixed value of $\mu$ is computed, possibly iterating over multiple primal-dual steps. Subsequently, the barrier parameter is updated and the computation continues by solution of the next barrier problem, starting from the approximate solution of the previous one. The approximate solution for the barrier problem \eqref{eq:Barrier}, for a given value of $\mu_j$, is required to satisfy the tolerance
\begin{equation}
    E_{\mu}(\bfx^{j+1}, \bfs^{j+1}, \lamE^{j+1}, \lamI^{j+1}, \bfz^{j+1}, \bfy^{j+1}) <  \kappa_{\epsilon} \mu_j
\end{equation}
for a constant $\kappa_{\epsilon}>0$ before the algorithm continues with the solution of the next barrier problem. The optimality error for the barrier problem is defined by considering the individual parts of the primal-dual equations \eqref{eq:OC}, that is, the dual feasibility (optimality), primal feasibility (constraint violations), and the complementarity conditions,
\begin{equation}
\label{eq:NLPerror}
  E_{\mu}(\bfx, \bfs, \lamE, \lamI, \bfz, \bfy) = 
  \max \left(
    \Vert\boldsymbol{l}_a\Vert_{\infty},
    \Vert\boldsymbol{l}_b\Vert_{\infty},
    \Vert\boldsymbol{l}_c\Vert_{\infty},
    \Vert\boldsymbol{l}_d\Vert_{\infty},
    \Vert\boldsymbol{l}_e\Vert_{\infty},
    \Vert\boldsymbol{l}_f\Vert_{\infty}
  \right).
\end{equation} 

In the monotone barrier update strategy, the new barrier parameter is obtained from
\begin{equation}
\label{eq:MUupdate}
    \mu_{j+1} = \max \left( \frac{\epsilon_{\text{tol}}}{10},\  \min \left(\kappa_{\mu} \mu_j,\  \mu_j^{\theta_{\mu}} \right) \right)
\end{equation}
with constants $\kappa_{\mu} \in (0,1)$ and $\theta_{\mu} \in (1,2)$. In this way, the barrier parameter is eventually decreased at a superlinear rate. On the other hand, the update rule \eqref{eq:MUupdate} does not allow $\mu$ to become smaller than necessary given the desired tolerance $\epsilon_{\text{tol}}$, thus avoiding numerical difficulties at the end of the optimization procedure. The monotone Fiacco--McCormick strategy can be very sensitive to the choice of the initial point, the initial value of the barrier parameter, and the scaling of the problem. Furthermore, different problems might favor strategies for selecting the barrier parameter at every iteration of an IP method, that is, for every primal-dual step computation.
Adaptive strategies commonly choose $\mu_{k+1}$ proportionally to the duality measure for the $k$th iterate,
\begin{equation}
    \mu_{k+1} = \sigma \tau_k,
\end{equation}
where $\sigma > 0$ is a centering parameter and $\tau$ denotes the duality measure \eqref{eq:dualMeasure}. The adaptive strategies vary in how the centering parameter is determined.  Two adaptive strategies implemented in \IPOPT{} are discussed next.

\subsubsection{Mehrotra's Predictor-Corrector \label{sec:adaptive1}}

Mehrotra's proposed  a predictor-corrector principle~\cite{muADAPTIVEmehrotra} for computing the search direction. The centering parameter is computed as the ratio between the duality measure \eqref{eq:dualMeasure} in the current iterate and the iterate updated by the predictor step, considering the longest possible step sizes that retain the nonnegativity of the variables in the barrier problem. If good progress in the duality measure is made in the predictor step, the centering parameter obtained in this way is small, $\sigma < 1$; therefore, the $\mu$ will be small in the next iteration. In other cases $\sigma$ may be chosen to be greater than $1$. This heuristic is based on experimentation with linear programming problems, and has proved to be effective for convex quadratic programming.

\subsubsection{Quality Function \label{sec:adaptive2}}
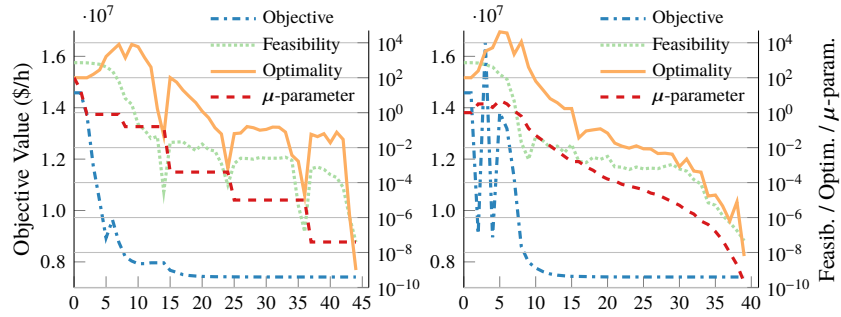
\begin{figure}[t]
\noindent
\definecolor{mycolor0}{HTML}{2B83BA}
\definecolor{mycolor1}{HTML}{D7191C}
\definecolor{mycolor2}{HTML}{FDAE61}
\definecolor{mycolor3}{HTML}{ABDDA4}
\begin{tikzpicture}

\pgfplotstableread{
iter objective infpr infdu mu
0 1.4580557e+07 7.23e+02 1.00e+02 100
1 1.4580557e+07 7.23e+02 9.93e+01 19.9526
2 1.3735969e+07 7.07e+02 9.89e+01 0.794328
3 1.1739171e+07 6.50e+02 1.56e+02 0.794328
4 1.0363599e+07 5.84e+02 2.82e+02 0.794328
5 8.9339278e+06 4.13e+02 1.64e+03 0.794328
6 9.6379213e+06 2.37e+02 3.66e+03 0.794328
7 8.7620973e+06 5.98e+01 8.07e+03 0.794328
8 8.2004783e+06 6.45e+00 1.38e+03 0.158489
9 8.0341158e+06 3.00e+00 8.01e+03 0.158489
10 7.9283934e+06 2.18e-01 5.90e+03 0.158489
11 7.9319774e+06 1.03e-01 1.42e+03 0.158489
12 7.9647282e+06 3.34e-02 4.03e+02 0.158489
13 7.9650992e+06 4.31e-02 1.68e+00 0.158489
14 7.9652231e+06 3.07e-05 4.80e-02 0.158489
15 7.6726267e+06 1.98e-02 1.01e+02 0.000398107
16 7.5928769e+06 2.08e-02 5.57e+01 0.000398107
17 7.5040049e+06 2.15e-02 1.92e+01 0.000398107
18 7.4694809e+06 1.14e-02 8.28e+00 0.000398107
19 7.4491674e+06 6.86e-03 3.72e+00 0.000398107
20 7.4347607e+06 1.58e-02 1.20e+00 0.000398107
21 7.4287828e+06 7.69e-03 5.11e-01 0.000398107
22 7.4228334e+06 6.39e-03 1.36e-01 0.000398107
23 7.4211611e+06 1.18e-03 6.02e-02 0.000398107
24 7.4210356e+06 7.30e-05 6.45e-04 0.000398107
25 7.4174057e+06 1.22e-03 6.32e-02 1e-05
26 7.4150901e+06 2.47e-03 6.81e-02 1e-05
27 7.4142527e+06 2.37e-03 1.60e-01 1e-05
28 7.4138881e+06 2.18e-03 1.37e-01 1e-05
29 7.4128300e+06 2.73e-03 9.91e-02 1e-05
30 7.4124917e+06 2.43e-03 1.10e-01 1e-05
31 7.4120861e+06 2.50e-03 1.49e-01 1e-05
32 7.4118893e+06 2.44e-03 1.46e-01 1e-05
33 7.4116955e+06 2.83e-03 7.73e-02 1e-05
34 7.4115888e+06 1.90e-03 3.70e-03 1e-05
35 7.4115969e+06 4.95e-06 1.66e-03 1e-05
36 7.4115962e+06 1.54e-07 1.87e-05 1e-05
37 7.4115153e+06 6.27e-04 5.78e-02 3.98107e-08
38 7.4114687e+06 7.77e-04 4.25e-02 3.98107e-08
39 7.4114470e+06 5.71e-04 6.23e-02 3.98107e-08
40 7.4114314e+06 2.84e-04 1.90e-02 3.98107e-08
41 7.4114296e+06 1.93e-04 7.56e-02 3.98107e-08
42 7.4114269e+06 5.34e-05 2.96e-02 3.98107e-08
43 7.4114261e+06 1.07e-06 1.37e-06 3.98107e-08
44 7.4114261e+06 3.92e-08 9.91e-10 3.98107e-08
}\datatable;


\begin{axis}[
width=0.47\linewidth, height=5cm,
ymin=0.7e7,
ymax=1.7e7,
xmin=0,
xmax=46,
axis lines*=left,
yticklabel style={
	/pgf/number format/fixed,
	/pgf/number format/fixed zerofill,
	/pgf/number format/precision=1
},
xtick={0,5,10,15,20,25,30,35,40,45,50,55,60,65,70,75,80,90},
yticklabel style={font=\scriptsize},
xticklabel style={font=\scriptsize},
ylabel style={font=\footnotesize},
xlabel style={font=\footnotesize},
ylabel={Objective Value (\$/h)},
y label style={at={(axis description cs:-0.11,.5)},anchor=south},
legend=none    
]
\addplot[color=mycolor0, very thick, dashdotted] table[x=iter, y=objective] {\datatable}; \label{plot:objE1k}
\end{axis}

\begin{axis}[name=symb,
width=0.47\linewidth, height=5cm,
ymode=log, 
log origin=infty,
ymin=1e-10,
ymax=5e4,
axis lines*=right, 
axis x line=none, 
ymajorgrids,
xmin=0,
xmax=46,
ytick={1e-10, 1e-8, 1e-6, 1e-4, 1.00e-2,  1e0,  1e2,  1e4},
xticklabel style={rotate=0, xshift=-0.0cm, anchor=north, font=\footnotesize},
yticklabel style={font=\scriptsize},
ylabel style={font=\footnotesize},
xlabel style={font=\footnotesize},
y label style={at={(axis description cs:1.30,.5)},anchor=south},
legend style={at={(1.0,1.12)},legend cell align=left,align=right,draw=none,fill=none,font=\scriptsize,legend columns=1}
]
\addlegendimage{/pgfplots/refstyle=plot:objE1k}\addlegendentry{Objective}

\addplot[color=mycolor3, very thick, densely dotted] table[x=iter, y=infpr] {\datatable};
\addlegendentry{Feasibility}

\addplot[color=mycolor2, very thick] table[x=iter, y=infdu] {\datatable};
\addlegendentry{Optimality}


\addplot[color=mycolor1, very thick,, dashed] table[x=iter, y=mu] {\datatable}; 
\addlegendentry{$\mu$-parameter}
  
\end{axis}

\end{tikzpicture}
 \hspace{-6mm} 
\definecolor{mycolor0}{HTML}{2B83BA}
\definecolor{mycolor1}{HTML}{D7191C}
\definecolor{mycolor2}{HTML}{FDAE61}
\definecolor{mycolor3}{HTML}{ABDDA4}
\begin{tikzpicture}

\pgfplotstableread{
iter objective infpr infdu mu
0 1.4580557e+07 7.23e+02 1.00e+02 1
1 1.4580539e+07 7.23e+02 9.93e+01 1
2 9.0895028e+06 6.87e+02 2.45e+02 3.16228
3 1.6527256e+07 5.51e+02 3.30e+03 3.16228
4 8.9605499e+06 4.13e+02 5.15e+03 1.99526
5 1.3813134e+07 1.54e+02 4.19e+04 5.01187
6 1.3188031e+07 8.29e+01 3.61e+04 3.16228
7 1.1266757e+07 4.71e+00 2.20e+03 1.25893
8 8.5970607e+06 1.96e-02 1.14e+04 0.630957
9 7.9902927e+06 2.12e-03 4.00e+02 0.125893
10 7.7648004e+06 4.37e-02 6.33e+01 0.0501187
11 7.6226252e+06 2.51e-02 2.03e+01 0.0251189
12 7.5321201e+06 9.95e-03 6.36e+00 0.0125893
13 7.4837337e+06 1.22e-02 3.72e+00 0.00630957
14 7.4574647e+06 5.29e-03 1.71e+00 0.00316228
15 7.4371022e+06 1.51e-02 1.72e+00 0.00158489
16 7.4342676e+06 1.63e-03 3.62e-02 0.00158489
17 7.4281029e+06 1.38e-03 8.96e-02 0.000794328
18 7.4225580e+06 2.80e-03 1.04e-01 0.000501187
19 7.4182674e+06 2.12e-03 1.19e-01 0.000251189
20 7.4156499e+06 3.14e-03 6.73e-02 0.000158489
21 7.4146304e+06 8.02e-04 1.81e-02 0.000125893
22 7.4137064e+06 7.61e-04 1.21e-02 7.94328e-05
23 7.4131317e+06 5.60e-04 9.34e-03 6.30957e-05
24 7.4125733e+06 7.72e-04 1.26e-02 5.01187e-05
25 7.4122727e+06 6.10e-04 8.54e-03 3.98107e-05
26 7.4120148e+06 6.73e-04 8.38e-03 2.51189e-05
27 7.4118589e+06 6.25e-04 4.58e-03 1.99526e-05
28 7.4117009e+06 8.59e-04 4.78e-03 1.25893e-05
29 7.4115799e+06 1.09e-03 4.22e-03 7.94328e-06
30 7.4115202e+06 5.96e-04 8.19e-04 5.01187e-06
31 7.4114712e+06 4.51e-04 2.22e-03 2.51189e-06
32 7.4114499e+06 1.47e-04 4.54e-04 1.58489e-06
33 7.4114347e+06 8.19e-05 3.96e-04 6.30957e-07
34 7.4114312e+06 7.00e-06 1.70e-05 3.98107e-07
35 7.4114279e+06 5.02e-06 1.93e-05 1.58489e-07
36 7.4114263e+06 1.32e-06 4.64e-06 5.01187e-08
37 7.4114257e+06 5.14e-07 6.31e-07 1e-08
38 7.4114255e+06 1.78e-07 7.91e-06 1.99526e-09
39 7.4114255e+06 4.73e-08 6.25e-09 1.99526e-10
}\datatable;


\begin{axis}[
width=0.47\linewidth, height=5cm,
ymin=0.7e7,
ymax=1.7e7,
xmin=0,
xmax=41,
axis lines*=left,
yticklabel style={
	/pgf/number format/fixed,
	/pgf/number format/fixed zerofill,
	/pgf/number format/precision=1
},
xtick={0,5,10,15,20,25,30,35,40,45,50,55,60,65,70,75,80,90},
yticklabel style={font=\scriptsize},
xticklabel style={font=\scriptsize},
ylabel style={font=\footnotesize},
xlabel style={font=\footnotesize},
y label style={at={(axis description cs:-0.11,.5)},anchor=south},
legend=none    
]
\addplot[color=mycolor0, very thick, dashdotted] table[x=iter, y=objective] {\datatable}; \label{plot:objE2k}
\end{axis}

\begin{axis}[name=symb,
width=0.47\linewidth, height=5cm,
ymode=log, 
log origin=infty,
ymin=1e-10,
ymax=5e4,
axis lines*=right, 
axis x line=none, 
ymajorgrids,
xmin=0,
xmax=41,
ytick={1e-10, 1e-8, 1e-6, 1e-4, 1.00e-2,  1e0,  1e2,  1e4},
xticklabel style={rotate=0, xshift=-0.0cm, anchor=north, font=\footnotesize},
yticklabel style={font=\scriptsize},
ylabel style={font=\footnotesize},
xlabel style={font=\footnotesize},
ylabel={Feasib. / Optim. / $\mu$-param. },
y label style={at={(axis description cs:1.30,.5)},anchor=south},
legend style={at={(1.0,1.12)},legend cell align=left,align=right,draw=none,fill=none,font=\scriptsize,legend columns=1}
]
\addlegendimage{/pgfplots/refstyle=plot:objE2k}\addlegendentry{Objective}

\addplot[color=mycolor3, very thick, densely dotted] table[x=iter, y=infpr] {\datatable};
\addlegendentry{Feasibility}

\addplot[color=mycolor2, very thick] table[x=iter, y=infdu] {\datatable};
\addlegendentry{Optimality}


\addplot[color=mycolor1, very thick, dashed] table[x=iter, y=mu] {\datatable}; 
\addlegendentry{$\mu$-parameter}
  
\end{axis}

\end{tikzpicture}
\caption{Barrier parameter update strategies (left: monotone $\mu_0=100$; right: adaptive).}
\label{fig:MUupdateStrategies}
\end{figure}

The adaptive barrier update strategy based on the quality function, as suggested in \cite{muADAPTIVEqualityfcn}, is trying to determine the centering parameter by minimizing a linear approximation of the quality function. The quality function is a measure defined by the infeasibility norms in the current iterate updated by the probing search direction, which is expressed as a function of the sought parameter $\sigma$. The minimization problem is solved by a golden bisection procedure on the specified $(\sigma_{\text{min}}, \sigma_{\text{max}})$ interval with a maximum of 12 bisections.
The evaluation of the barrier update strategies on both linear and nonlinear problems revealed superior performance of the adaptive methods over the monotone strategy, both in terms of CPU time and number of IP iterations. Although the results were more pronounced on the linear benchmarks, significant improvements can be expected by using adaptive strategies, particularly in applications where the function evaluation has the dominant cost \cite{muADAPTIVEqualityfcn}. Figure \ref{fig:MUupdateStrategies} depicts the convergence with different barrier parameter update strategies. The value of the barrier parameter $\mu$ over the iterations of the IP is shown for the two update strategies. Feasibility, optimality and the objective function are shown as well. The convergence tolerance for both benchmarks was set to $\text{tol}=0.01$.


\subsection{Problem Scaling and Convergence Criteria  \label{sec:Scaling}}
Optimal control of realistic industrial and engineering problems, such as modern power networks, multienergy carrier systems, the variables and constraints encountered, commonly involve different scales that usually differ by several orders of magnitude. Sophisticated scaling is necessary to 
remedy problems related to establishing accurate stopping criteria, improving convergence deteriorated by unbalanced direction vectors, and dealing with loss of accuracy of the descent direction computation due to bad conditioning of the associated KKT systems. In the ideal case, not only the variables but also the functions should be scaled so that changing a variable by a given amount has a comparable effect on any function which depends on these variables or, in other words, so that the nonzero elements of the function gradients are of the same order of magnitude. For this purpose, gradient-based scaling is commonly employed so that at the starting point the gradients are scaled close to one. The scaling factors for the gradients are defined as
\begin{align}
    s_f   &= \min (1,\ g_{\text{max}} / \Vert\nabla_x f(\bfx_0)\Vert_{\infty}), \\
    s_g^{(j)} &= \min (1,\ g_{\text{max}} / \Vert\nabla_x \ce^{(j)}(\bfx_0)\Vert_{\infty}),\  j=1\ldots \ne, \\
    s_h^{(j)} &= \min (1,\ g_{\text{max}} / \Vert\nabla_x \ci^{(j)}(\bfx_0)\Vert_{\infty}),\ j=1\ldots \ni,
\end{align}
for a given $g_{\text{max}} > 0$. If the maximum gradient is above this value, then gradient-based scaling will be performed. Note that all gradient components in the scaled problem are at most of size $g_{\text{max}}$ at the starting point. The scaling factors are computed only at the beginning of the optimization using the starting point and kept constant throughout the whole optimization process.

Even if the original problem is well scaled, the multipliers $\lamE, \lamI, \bfz$ might become very large, for example, when the gradients of the active constraints are (nearly) linearly dependent at a solution of \eqref{eq:NLP}. In this case, the algorithm might encounter numerical difficulties satisfying the unscaled primal-dual equations \eqref{eq:SCOPFhessianSymmetric} to a tight tolerance. The convergence criteria in \eqref{eq:NLPerror}, therefore, need to be scaled accordingly. The scaled optimality error used to determine the convergence criteria is defined as
\begin{equation}
\label{eq:NLPerrorScaled}
  E_0(\bfx, \bfs, \lamE, \lamI, \bfz) = 
  \max \left(
    \frac{\Vert\boldsymbol{l}_a\Vert_{\infty}}{s_1},
    \frac{\Vert\boldsymbol{l}_b\Vert_{\infty}}{s_1},
    \Vert\boldsymbol{l}_c\Vert_{\infty},
    \Vert\boldsymbol{l}_d\Vert_{\infty},
    \frac{\Vert\boldsymbol{l}_e\Vert_{\infty}}{s_2},
    \frac{\Vert\boldsymbol{l}_f\Vert_{\infty}}{s_2}
  \right),
\end{equation}
where the scaling factors $s_1, s_2$ are defined as
\begin{equation}
    s_1 = \dfrac{\max \left( s_{\text{max}}, \frac{\Vert\lamE\Vert_1 + \Vert\lamI\Vert_1 + \Vert\bfz\Vert_1 + \Vert\bfy\Vert_1}{\ne + \ni + N_x + \ni} \right)}{s_{\text{max}}}, \quad s_2 = \dfrac{\max \left( s_{\text{max}}, \frac{\Vert\bfz\Vert_1+\Vert\bfy\Vert_1}{N_x+\ni} \right)}{s_{\text{max}}}.
\end{equation}

The overall \IPOPT{} algorithm terminates successfully, if the NLP error for the current iterate with $\mu = 0$ in \eqref{eq:NLPerrorScaled},
\begin{equation}
  \label{eq:NLPerrorZERO}
  E_{0}(\bfx, \bfs, \lamE, \lamI, \bfz, \bfy) \leq \epsilon_{\text{tol}},
\end{equation}
becomes smaller than the user provided value $\epsilon_{\text{tol}} >0$, and if the individual criteria according to dual, primal, and complementarity conditions in \eqref{eq:NLPerrorScaled} are met. Each criterion uses a separate, user provided tolerance value.
\section{IP Methods for OPF Problems}

Recent developments in modern power grids involve widespread deployment of intermittent renewable generation, embrace installation of a wide variety of energy storage devices, as well as an increasing and widespread usage of electric vehicles. These developments will motivate fundamental changes in methods and tools for the optimal daily operation and planning of modern power grids. Operational decisions taken by power system operators on a daily basis are commonly assisted by repeatedly solving OPF problems, aiming to determine optimal operating levels for electric power plants, so that the overall electricity generation cost is minimized, while at the same time it satisfies load demands imposed throughout the transmission grid and meets safe operating limits. In actual industrial operations the entire distribution network needs to be optimized in real time, approximately every five minutes according to several independent system operators to ensure variations in load demand, renewable generation, and real-time electricity market responses to electricity prices are accurately met. 

\subsection{Optimal Power Flow}
The OPF problem seeks a solution that minimizes the cost of the electricity generation $f$, while satisfying the power flow balance, maximum power flow over the transmission lines, and the bounds of the bus voltages and the generator limits.
Consider a power network with $\Nb$ buses, $\Ng$ generators, and $\Nl$ transmission lines. The bus voltage vector $\vc \in \ComplexSet{\Nb}$ is defined in polar notation as $\vc = \bfv e ^{j\btheta}$, where  $\bfv, \btheta \in \R{\Nb}$ specify the magnitude and phase of the complex voltage. The complex voltages $\vc$ determine the entire network power flow that can be computed using the Kirchhoff equations, network configuration, and properties of its components. The magnitude of the voltage components is bounded by the limits \eqref{OPF:VoltageLimits}, while the phase is determined relative to a single reference bus.
The current injections, $\ic \in \ComplexSet{\Nb}$, are defined as $\ic = \Ybc \vc$, where $\Ybc \in \ComplexSet{\Nb \times \Nb}$ is the bus admittance matrix. The complex power at each bus of the network, $\sc = \vc \ic^*$, $\sc \in \ComplexSet{\Nb}$, and the power demand consumption $\sdc \in \ComplexSet{\Nb}$, are to be balanced by the net power injections from the generators $\sgc \in \ComplexSet{\Ng}$. Thus, the AC nodal power flow balance equations \eqref{OPF:PowerFlow} are expressed as a function of the complex bus voltages and generator injections as $\PowerFlow := \sc + \sdc - \Cgb\sgc= \0$, where $\Cgb \in \R{\Nb \times \Ng}$ is the generator connectivity matrix.

Generator power injections $\sgc = \bfp + j \bfq$ are expressed in terms of real and reactive power components $\bfp, \bfq\in \R{\Ng}$, respectively. The output of the generators is limited by the lower and upper bounds \eqref{OPF:PowerLimitsP} and \eqref{OPF:PowerLimitsQ}. Each bus has an associated complex power demand $\sdc$, which is assumed to be known at all of the buses and is modeled by a static polynomial (ZIP) model~\cite{matpowerManual}. If there are no loads connected to the bus $i$ then $\{\sdc\}_i = \0$.
Real-world transmission lines are limited by the instantaneous amount of power that can flow through the lines due to the thermal limits \eqref{OPF:LineFlowLimits}. The apparent power flow in the transmission lines, $\sfc \in \ComplexSet{\Nl}$ and $\stc \in \ComplexSet{\Nl}$, are therefore limited by the power injections at both ends of the lines, which cannot exceed a prescribed upper bound $\Max{\F}_\text{L}$. 
The ``from" and ``to" ends of the line, denoted as $f$ and $t$, respectively, specify the buses that are connected to the corresponding ends of the line. Squared values of the apparent power magnitude are usually used in practice, such that $\LineFlow := \sfc(\sfc)^* \leq (\Max{\F}_\text{L})^2$. 
Overall, the OPF problem is formulated as
\begin{linenomath*}
	\begin{subequations}
		\label{OPF}
		\begin{align}
			\minimize{\btheta,\bfv,\bfp,\bfq} \;  & \sum_{l=1}^{\Ng} f_l(\pes{l}{}) \label{OPF:objective}\\ 
			\subject{} \; & \PowerFlow(\btheta, \bfv, \bfp, \bfq) = \0, \label{OPF:PowerFlow}\\  
			& \LineFlow(\btheta, \bfv)  \leq \Max{\F}_\text{L}, 	\label{OPF:LineFlowLimits}\\
			& \Min{\bfv}      \leq \bfv      \leq \Max{\bfv},\quad \btheta^{ref} = 0,  \label{OPF:VoltageLimits}\\
			& \Min{\bfp} \leq \bfp  \leq \Max{\bfp}, \label{OPF:PowerLimitsP}\\ 
			& \Min{\bfq} \leq \bfq  \leq \Max{\bfq}. \label{OPF:PowerLimitsQ}
		\end{align}
	\end{subequations}
\end{linenomath*}
The presented AC steady-state power grid model is following \MATPOWER{}~\cite{matpower}.

\subsection{Structure-Exploiting IP Methods -- Security Constrained and Multiperiod OPF}
Real-world real-time implementation of OPF problems for energy systems, still remain computationally intractable. This is mainly for two reasons. The real world OPF problem is time coupled, owing to the presence of smart loads and energy storage devices such as batteries for demand shaping and deferral. Additional time couplings of the OPF problem at each time period are introduced by generator ramp rate limits. The higher the number of time periods considered, the larger the resulting optimal control problem becomes. For a significantly large number of time periods (each of five-minutes length) the problem becomes notoriously difficult to solve and for this purpose several approximations and simplifications are currently employed by the industry in order to meet real-time responses.
Furthermore, the system operators have to foresee possible contingency events and operate the grid in a such a way that its operation will remain secure in the event of any contingencies.

Grid security is the focus of the SCOPF problem~\cite{SCOPFinvention, SCOPFcorrective}, which seeks an optimal solution that  remains feasible under any postulated contingency event,  thus making the grid operation secure.  It supplements the standard OPF problem with constraints for the nodal power flow balance  \eqref{SCOPF:PowerFlow}, the branch flow limits \eqref{SCOPF:LineFlowLimits}, and other operational limits \eqref{SCOPF:VoltageLimits},  \eqref{SCOPF:PowerLimits}, which have to be honored not only for the nominal case $c_0$, but also for every contingency event $c \in \mathcal{C}$, $N_c = |\mathcal{C}|$, such as a generator or a transmission line failure. An increase of the number of considered contingencies requires the introduction of additional variables and constraints that in turn result in a significant problem size growth, rendering it computationally intractable for standard general purpose optimization tools. 
The contingencies are modeled by the admittance matrices $\Ybc_c$, which are updated accordingly for each scenario.
The values of the control variables are coupled in all system scenarios, as expressed by the two nonanticipatory constraints \eqref{SCOPF:SCOPFVm} and \eqref{SCOPF:SCOPFPg}. These declare that the voltage magnitude and real power generation at the PV buses $\mathcal{B}_{PV}$ should remain the same as in the nominal scenario $c_0$, regardless of which contingency they are associated with. The only generator that is allowed to change its output is the generator at the singleton reference bus $\mathcal{B}_{ref}$, as its real power generation can be modified to refill the power transmission losses occurring in each contingency $c$. 

\begin{figure}[t]
\begin{minipage}[t]{0.495\textwidth}
\begin{linenomath*}
\begingroup\leqnos 
\begin{subequations}
	\label{SCOPF}
	\begin{align}
	 \nonumber 
	\minimize{\btheta_{c},\bfv_{c},\bfp_{c},\bfq_{c}} \; & \sum_{l=1}^{\Ng} f_l(\pes{l}{0}) \\ 
	\subject{} \; & \forall c \in \{ c_0, c_1, \ldots, c_{N_c}\}, \nonumber \\
	\label{SCOPF:PowerFlow}
	& \PowerFlowSecurity(\btheta_c, \bfv_c, \bfp_c, \bfq_{c}) = \0, \\  
	\label{SCOPF:LineFlowLimits}
	& \LineFlowSecurity(\btheta_c, \bfv_c)  \leq \Max{\F}_\text{L}, \\
	\label{SCOPF:VoltageLimits}
	& \Min{\bfv}      \leq \bfv_c      \leq \Max{\bfv},\\
	& \btheta_c^{ref} = 0,  \\
	\label{SCOPF:PowerLimits}
	& \Min{\bfp} \leq \bfp_c  \leq \Max{\bfp}, \\ 
	& \Min{\bfq} \leq \bfq_c  \leq \Max{\bfq}, \\
	\label{SCOPF:SCOPFVm}
	& \forall b \in \mathcal{B}_{PV}: \bfv_c = \bfv_{c_0},\\
	\label{SCOPF:SCOPFPg}
	& \forall  g \in \mathcal{B}_{PV}: \bfp_c = \bfp_{c_0}.
	\end{align}
\end{subequations} \endgroup
\end{linenomath*}
\end{minipage}
\begin{minipage}[t]{0.50\textwidth}
\begin{linenomath*}
\begin{subequations}
	\label{MPOPF}
	\begin{align}
    \nonumber 
	\minimize{\btheta_{n},\bfv_{n},\bfp_{n},\bfq_{n}} \; & \sum_{l=1}^{\Ng} f_l(\pes{l}{0}) \\ 
	\subject{}  \; & \forall n \in \{ 1, 2, \ldots, N\}, \nonumber \\
	\label{MPOPF:PowerFlow}
	& \PowerFlowSecurity(\btheta_n, \bfv_n, \bfp_n, \bfq_{n}) = \0, \\  
	\label{MPOPF:LineFlowLimits}
	& \LineFlowSecurity(\btheta_n, \bfv_n)  \leq \Max{\F}_\text{L}, \\
	\label{MPOPF:VoltageLimits}
	& \Min{\bfv}      \leq \bfv_n      \leq \Max{\bfv},\\
	& \btheta_n^{ref} = 0,  \\
	\label{MPOPF:PowerLimits}
	& \Min{\bfp} \leq \bfp_n  \leq \Max{\bfp}, \\ 
	\label{MPOPF:MPOPFVm}
	& \Min{\bfq} \leq \bfq_n  \leq \Max{\bfq}, \\
	\label{MPOPF:energy}
	& \bepsilon^{\textbf{min}} \leq \bepsilon_n \leq \bepsilon^{\textbf{max}}.
	\end{align}
\end{subequations}
\end{linenomath*}
\end{minipage}
\caption{SCOPF (left) and MPOPF (right) problem formulations.}
\end{figure}

Time-coupled formulations, such as storage scheduling, or storage placement, are collectively known as MPOPF problems \eqref{MPOPF}. Similar to the SCOPF, addition of a large number of time periods results in problem size growth, rendering it computationally intractable~\cite{beltistos}. The OPF constraints must hold in each time period, and the inter-temporal coupling is introduced by energy storage devices and  generator ramp limits. For a practical MPOPF application, consider $\Ns$ energy storage units. Each storage unit in the network is modeled by two network power injections for each time period $n$.
A positive active power injection $\psdi_n \in \R{}$, $\psdi_n \geq 0$  models the discharging of storage unit $i$. A negative active power injection $\psci_{n} \in \R{}$, $\psci_n \leq 0$ models the charging of storage unit $i$. The vector of active storage power injections $\ps_n \in \R{2\Ns}$ is defined as 
\begin{equation}
\ps_n = (\psdni{1}{n}, \cdots ,\psdni{\Ns}{n},  \pscni{1}{n},  \cdots , \pscni{\Ns}{n})
\end{equation}
and bounded by $\psmin \leq \ps_n \leq \psmax$. Identical definitions apply for the reactive storage power injections $\qsdni{i}{n}$, $\qscni{i}{n}$,  $\qs_n$ with bounds  $\qsmin$ and $\qsmax$ .
Together, they yield the complex storage power injections $\sS_n = \ps_n + j\qs_n$. Similarly, $\sg_n = \pg_n +j\qg_n$ is a vector of generator power injections.
The complex power at each bus must be balanced by the power demand $\sdc_n$  and the vector of free complex power injections  
\begin{equation}
\scomplex_n = \pmat{\sg_n \\ \sS_n} = 
\underbrace{\pmat{\pg_n \\ \ps_n}}_{\bfp_n}  + j \underbrace{\pmat{\qg_n \\ \qs_n}}_{\bfq_n}
\end{equation}
in each time period, as specified by the constraint \eqref{MPOPF:PowerFlow}.
The evolution of the vector of storage levels $\bepsilon_n \in\R{\Ns}$ follows the update equation
\begin{align}
	\label{Eq:linup}
	\bepsilon_{n+1} = \bepsilon_n + \Bs \; \ps_n \quad n = 0,1,\ldots,N-1,
\end{align}
and introduces a coupling between the individual time periods. The energy level in each period needs to honor the storage capacity, as expressed by the constraint \eqref{MPOPF:energy}.
The initial storage level is denoted $\bepsilon_0$ and the constant matrix $\Bs \in\R{\Ns\times2\Ns}$ models discharging and charging efficiencies of the storage devices.
\subsection{Impact of Slack Variables Elimination}
\begin{figure}[t]
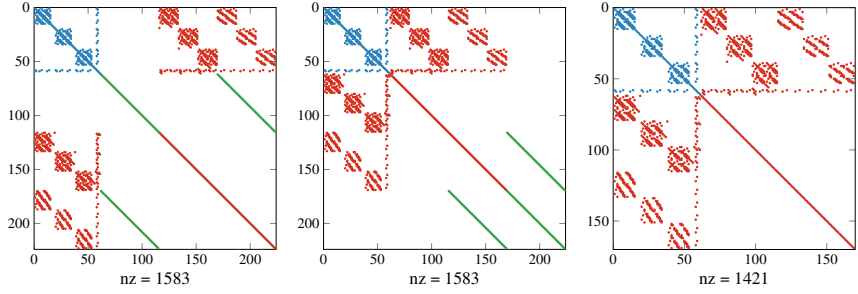

\begin{center}
  \input{SparseMatrices/KKT_SCOPFcase9ns3.tex}
  ~
  \input{SparseMatrices/KKT_SCOPFcase9ns3_slack.tex}
  ~ 
  \input{SparseMatrices/KKT_SCOPFcase9ns3_Reduced.tex}
\caption{Structure of the SCOPF KKT system \eqref{eq:SCOPFhessianSymmetric} with two contingencies, reordered according to \eqref{eq:SCOPFslackAtEnd}, and, finally, the reduced KKT with the slacks removed \eqref{eq:SCOPFhessianReduced}.}
\label{fig:KKTsystemSCOPF}
\end{center}
\end{figure}

Figure~\ref{fig:KKTsystemSCOPF} illustrates the symmetric KKT structure of the SCOPF problem for a simple power grid, together with the reduced variant, where the slack variables are eliminated. Realistic power grids are significantly larger and contain proportionally more nonzero entries, but the structure remains very similar.
The expected benefits of solving the reduced KKT system compared to the original system are savings both in terms of memory requirements for storing the sparse $\L$ factor of the $\L D\L^\transpose$ factorization of the symmetric indefinite system, and possibly faster factorization and solution times due to a smaller number of required floating point operations. The numerical evaluation of the benefits of solving the reduced system are summarized in Figure~\ref{fig:slackEliminationResults}.
The elimination of the slack variables from the KKT system reduces its dimension by approximately $30$\% with $13$\% fewer nonzeros in the KKT system and up to $12$\% fewer nonzeros in the $\L$ factor, resulting in up to $28$\% memory savings, with similar reduction in solution time. Since in the neighborhood of the optimal solution some of the diagonal terms in $L_s$ approach zero, the associated slacks variables whose coefficients in $L_s$ are close to machine epsilon are not eliminated and are left to be treated by the direct sparse solver. This prevents the excessive ill-conditioning of the reduced system.

\section{Structure-Exploiting Solution Strategies for IP Optimization\label{chapter:solution} \label{sec:schur}}

\begin{figure}[t]
\begin{minipage}[t]{0.465\textwidth}
	\begin{tikzpicture}

\pgfplotstableread{
id	case	order	memory	time
1	PEGASE1k	metis	28.11	21.91
2	2737sop	metis	26.88	20.56
3	PEGASE2k	metis	28.56	16.73
4	3375wp	metis	26.65	20.91
5	6495rte	metis	28.13	24.06
6	6515rte	metis	28.2	25.13
7	PEGASE9k	metis	19.2	22.88
8	ACTIVSg10k	metis	13.13	22.64
}\datatable;

\begin{axis}[name=symb,
width=\columnwidth, height=4cm, 
ybar,
bar width=3pt, 
ymin=0,
ymax=32,
xmin=0.5,
xmax=8.5,
axis lines*=left, 
ymajorgrids, yminorgrids,
xticklabels from table={\datatable}{case},
xtick={1, 2, 3, 4, 5, 6, 7, 8},
xticklabel style={align=center, rotate=45, xshift=-0.4cm, anchor=north, font=\scriptsize},
ytick={5,10,15,20,25,30},
yticklabel style={font=\scriptsize},
ylabel style={font=\scriptsize},
xlabel style={font=\scriptsize},
ylabel={Improvement (\%)},
legend style={at={(1.1,1.25)},fill=white,legend cell align=left,align=right,draw=none,font=\scriptsize,legend columns=4}
]
\addplot[fill=mycolor0] table[x=id, y=memory] {\datatable};
\addlegendentry{Memory requirements}
\addplot[fill=mycolor1, postaction={pattern=north east lines}] table[x=id, y=time] {\datatable}; 
\addlegendentry{Solution time}
\end{axis}

\end{tikzpicture}
	\caption{Improvement rate considering elimination of the slack variables. \label{fig:slackEliminationResults}}
\end{minipage}
\quad
\begin{minipage}[t]{0.48\textwidth}
\begin{tikzpicture}

\pgfplotstableread{
	ID  NPROC   MPTTOTAL    TTOTAL  TSOLVE  TFEVAL  MISC    CASE
	1   1  2.63E+05 2.62E+05    2.60E+05    9.99E+02    7.81E+02 PEGASE\\1354-2048
	2   1   2.63E+05    2.62E+05    2.60E+05    9.99E+02    7.81E+02    PEGASE\\9241-512
	3   1   3.40E+04   81519.712	81299.408	143.552	76.752  PEGASE\\9241-256
	4   1   4.73E+04    47276.652   46864.292   339.184 73.176  PEGASE\\13659-128
	5   1   7.40E+03    7369.956    7165.596    169.636 34.724  PEGASE\\13659-64
}\datatable;

\begin{axis}[name=symb,
width=\columnwidth, height=4cm, 
ybar,
bar width=3pt, 
ymode=log, 
log origin=infty,
ymin=10,
ymax=1e6,
axis lines*=left, 
ymajorgrids, yminorgrids,
xticklabels from table={\datatable}{CASE},
xtick={1, 2, 3, 4, 5, 6, 7, 8},
xticklabel style={align=center, rotate=45, xshift=-0.5cm, anchor=north, font=\scriptsize},
ytick={10,100, 1000, 10000, 100000, 1e6},
yticklabel style={font=\scriptsize},
ylabel style={font=\scriptsize},
xlabel style={font=\scriptsize},
ylabel={Time (s)},
legend style={at={(1.0,1.3)},fill=white,legend cell align=left,align=right,draw=none,font=\scriptsize,legend columns=2}
]
\addplot[fill=mycolor0] table[x=ID, y=TTOTAL] {\datatable};
\addlegendentry{Overall time}
\addplot[fill=mycolor1, postaction={pattern=north east lines}] table[x=ID, y=TSOLVE] {\datatable}; 
\addlegendentry{KKT solution}
\addplot[fill=mycolor2, postaction={pattern=north west lines}] table[x=ID, y=TFEVAL] {\datatable}; 
\addlegendentry{Function eval.}
\addplot[fill=mycolor3, postaction={pattern=crosshatch}] table[x=ID, y=MISC] {\datatable}; 
\addlegendentry{Other}
\end{axis}

\end{tikzpicture}
    \caption{Computational complexity of the IP method components.}
    \label{fig:complexity}
\end{minipage}
\end{figure}

Computers have evolved significantly over the past decade, at an even faster pace than modern power grids. Multicore and many-core computer architectures and distributed compute clusters are ubiquitous today, while at the same time no significant performance gains are expected for sequential codes due to faster clock frequencies of modern processors. Significant performance gains, however, may be achieved by algorithmic redesign tailored to the particular application that is also able to utilize multicore and many-core architectures with deep memory hierarchies. More importantly, the practical efficiency of the IP algorithms  highly depends on the linear algebra kernels used.
For large-scale optimal control problems, the computation of the search direction \eref{eq:SCOPFhessianSymmetric} determines the overall runtime. Hence, any attempt at accelerating the solution should be focused on the efficient solution of the KKT linear system.
In Figure~\ref{fig:complexity} we demonstrate how various IP method components contribute to 
the overall time for various OPF benchmarks. 
The number of IP iterations was fixed to five. Note that the solution of the linear system represents the majority of the overall time.

\subsection{Revealing the Structure of SCOPF and MPOPF Problems}
A widespread approach for solving KKT systems consists of employing black-box techniques such as  direct sparse solvers, due to their accuracy and robustness. The direct sparse solvers obtain the solution of the linear system by factorization and subsequent forward-backward substitutions. The factorization is a computationally expensive operation commonly introducing significant fill-in, which may quickly exhaust available memory on shared memory machines for large-scale linear systems. Furthermore, these solvers are not aware of the underlying structural properties of the KKT systems arising from many engineering problems which make it possible to significantly decrease time to solution by employing structure-exploiting algorithms and distributed memory computers. 

The appropriate structure emerges from the fact that each of the variables in the SCOPF optimization vector
$( \bfx, \lamE, \lamI)$ 
or the MPOPF optimization vector $( \bfx, \lamE, \lamI, \lamA)$
correspond to some contingency scenario $c = 0, 1, \ldots, N_c$, 
or the time period $n = 1, 2, \ldots, N$:

\noindent
\begin{minipage}{0.49\textwidth}
\vspace{-0.4cm}
\begingroup\leqnos 
\begin{align}
	\bfx &= (\bfx_0,\ldots,\bfx_{N_c},\bfx_g), \\ 
	\lamE &= ({\lamE}_0, \ldots, {\lamE}_{N_c}), \\
	\lamI &= ({\lamI}_0, \ldots, {\lamI}_{N_c}),
\end{align} \endgroup
\end{minipage}
\begin{minipage}{0.49\textwidth}
\vspace{-0.4cm}
\begin{align}
	\bfx &= (\bfx_0,\ldots,\bfx_{N}), \\ 
	\lamE &= ({\lamE}_0, \ldots, {\lamE}_{N}), \\
	\lamI &= ({\lamI}_0, \ldots, {\lamI}_{N}).
\end{align}
\end{minipage}

\medskip
\noindent
In order to reveal the scenario-local structure of the Hessian~\eqref{eq:SCOPFhessianReduced}, the variables corresponding to the same contingency are grouped together, i.e.,

\smallskip
\noindent
\begin{minipage}{0.49\textwidth}
\begingroup\leqnos 
\begin{equation}
	\bfu_c = (\bfx_c,\ {\lamE}_c,\ {\lamI}_c), 
\end{equation} \endgroup
\end{minipage}
\begin{minipage}{0.49\textwidth}
\begin{equation}
	\bfu_n = (\bfx_n,\ {\lamE}_n,\ {\lamI}_n), 
\end{equation}
\end{minipage}

\medskip
\noindent
and, thus, the global ordering will be

\smallskip
\noindent
\begin{minipage}{0.49\textwidth}
\begingroup\leqnos 
\begin{equation}
\bfu = (\bfu_0, \ldots, \bfu_{N_c}, \bfu_g), \label{eq:orderingSCOPF}
\end{equation} \endgroup
\end{minipage}
\begin{minipage}{0.49\textwidth}
\begin{equation}
\bfu = (\bfu_0, \ldots, \bfu_{N}, \bfu_g), \label{eq:orderingMPOPF}
\end{equation}
\end{minipage}

\medskip
\noindent 
where the coupling variables $\bfu_g$ are placed at the end of the new optimization vector $\bfu$. Coupling in the SCOPF problem, $\bfu_g = \bfx_g$, is introduced by the two nonanticipatory constraints \eqref{SCOPF:SCOPFVm} and \eqref{SCOPF:SCOPFPg}. The coupling in a case of the MPOPF problem, $\bfu_g = \lamA$, is introduced by the linear energy constraints \eqref{MPOPF:energy}.
Under the new orderings~\eqref{eq:orderingSCOPF} and~\eqref{eq:orderingMPOPF}, the Hessian matrix of the system~\eqref{eq:SCOPFhessianReduced} obtains the arrowhead structure (also described as bordered block-diagonal~\cite{Duff} or  dual block-angular~\cite{pardiso}) structure, as illustrated in Figures \ref{fig:blockPermutaitionSCOPF} and \ref{fig:blockPermutaitionMPOPF},
         
\begin{equation}
\label{eq:hessianPerm}
\begin{pmatrix}
	\A_0 &        &            &                & \B_0^\transpose \\
		   & \A_1 &	           &                & \B_1^\transpose \\
		   &        & \ddots &                & \vdots \\
		   & 	    &            &   \A_{N_c}     & \B_{N_c}^\transpose\\
	\B_0  & \B_1 &   \ldots &   \B_{N_c}   &   \C
\end{pmatrix}
\begin{pmatrix}
	\Delta \bfu_0\\
	\Delta \bfu_1 \\
	\vdots\\
	\Delta \bfu_n\\
	\Delta \bfu_g
\end{pmatrix} = 
\begin{pmatrix}
	\b_0\\
	\b_1 \\
	\vdots\\
	\b_n\\
	\b_C
\end{pmatrix},
\end{equation}
where the block matrices $\A_i$,

\noindent
\begin{minipage}{0.45\textwidth}
\begingroup\leqnos 
\begin{equation}
	\A_i = \begin{pmatrix}
	\tilde{\Hess}_{x_i,x_i} & \JeXX{_i,x_i}{\transpose} & \JiXX{_i,x_i}{\transpose}  \\
	\JeXX{_i,x_i}{} & \0 & \0 \\
	\JiXX{_i,x_i}{} & \0  & -L_{s_i}^{-1}
	\end{pmatrix},
\end{equation} \endgroup
\end{minipage}
\begin{minipage}{0.54\textwidth}
\begin{equation}
	\A_i = \begin{pmatrix}
	\tilde{\Hess}_{x_i,x_i} & \JeXX{_i,x_i}{\transpose} & \JiXX{_i,x_i}{\transpose}  & \0\\
	\JeXX{_i,x_i}{} & \0 & \0 & \0\\
	\JiXX{_i,x_i}{} & \0  & -L_{s_i}^{-1} & \0 \\
	\0 & \0 & \0 & \L_{A_i}
	\end{pmatrix},
\end{equation}
\end{minipage}

\smallskip
\noindent
incorporate the Hessian of the Lagrangian with respect to the scenario-local variables $\tilde{\Hess}_{x_i,x_i} = \nabla^2_{x_i x_i}\mathcal{L} + {X}_i^{-1} {Z}_i$  and the Jacobians of the constraints for the $i$th scenario with respect to the local variables $\JeXX{_i,x_i}{} = \nabla_{x_i} {\PowerFlow}_i$ and $\JiXX{_i,x_i}{\transpose} = \nabla_{x_i} {\LineFlow}_i$, as well as the diagonal entries corresponding to the eliminated slack variables. In the case of the SCOPF problem, the block $\C = \nabla^2_{x_g x_g}\mathcal{L} + X_g^{-1} Z_g$ contains Hessian of the Lagrangian with respect to the coupling variables $\bfx_g$, while in the case of the MPOPF problem it is a block of zeros. The off-diagonal blocks in the arrowhead SCOPF system are
\begin{equation}
	\B_i = \begin{pmatrix}
	\tilde{\Hess}_{x_g,x_i} \\
	\JeXX{_i,x_g}{\transpose}  \\
	\JiXX{_i,x_g}{\transpose}
	\end{pmatrix}^\transpose, \quad
	\B_i^\transpose = \begin{pmatrix}
	\tilde{\Hess}_{x_i,x_g} \\
	\JeXX{_i,x_g}{} \\
	\JiXX{_i,x_g}{}
	\end{pmatrix},
\end{equation}
where $\tilde{\Hess}_{x_i,x_g} = \nabla^2_{x_ix_g}\mathcal{L}$ represents the off-diagonal blocks of the Hessian of Lagrangian with respect to the local and coupling variables and $\JeXX{_i,x_g}{} = \nabla_{x_g} {\PowerFlow}_i$ and $\JiXX{_i,x_g}{} = \nabla_{x_g} {\LineFlow}_i$ are the Jacobians of the $i$th scenario with respect to the coupling variables.
The MPOPF coupling matrices $\B_1,\B_2,\ldots,\B_N \in \R{N \Ns \times N_A}$, where $N_A$ is the size of the diagonal blocks in \eqref{eq:hessianPerm}, contain the constant subblocks, which arise from the particular form of the linear constraints \eqref{MPOPF:energy}
\begin{align} \label{Eq:B12N} 
	\B_1 = \pmat{ \C_1                          \\ 
		\C_0                                           \\ 
		\C_0                                            \\ 
		\vdots                                          \\ 
		\C_0             },
	\B_2 = \pmat{ \0                             \\ 
		\C_1                                            \\ 
		\C_0                                            \\ 
		\vdots                                         \\ 
		\C_0             }, \ldots,        
	\B_N =&\pmat{ \0                          \\
		\0                                              \\ 
		\vdots                                        \\ 
		\0                                                \\ 
		\C_1             }.
\end{align} 

\begin{figure}[t]
     \centering
     \input{SparseMatrices/KKT_SCOPFcase9ns3_blocksReduced.tex}
    \qquad \qquad
   \input{SparseMatrices/KKT_SCOPFcase9ns3_arrowReduced.tex}
    \caption{Symmetrized SCOPF system~\eqref{eq:SCOPFhessianReduced} permuted to the arrowhead structure~\eqref{eq:hessianPerm}.}
    \label{fig:blockPermutaitionSCOPF}
     \centering
     \input{SparseMatrices/KKT_MPOPFcase9nt3ns2_blocksReduced.tex}
    \qquad \qquad
   \input{SparseMatrices/KKT_MPOPFcase9nt3ns2_arrowReduced.tex}
    \caption{Symmetrized MPOPF system~\eqref{eq:SCOPFhessianReduced} permuted to the arrowhead structure~\eqref{eq:hessianPerm}.}
    \label{fig:blockPermutaitionMPOPF}
\end{figure}

\subsection{Schur Complement Decomposition}

The direct factorization of the full KKT system is not feasible for large-scale SCOPF problems due to their growing size with the number of contingencies and associated factorization fill-in that quickly exhausts the available memory. Instead, the solution is obtained by a sequence of partial block elimination steps, which are decoupled, aiming to form the Schur complement of the system. This way, we detour the factorization of  the full KKT system, by factorizing only the smaller diagonal blocks as described in the Algorithm~\ref{alg:algo:bw}.
At the first step, the Schur complement $\bfS$ is formed,
\begin{equation}
	\bfS = \C - \sum_{i=0}^{N_c} \B_i \A_i^{-1} \B_i^\transpose, \label{eq:SCform}
\end{equation}
which in the general case becomes a dense matrix. Because the size of the coupling stays constant, independently of the number of contingency scenarios, the size of the Schur complement does not increase with an increasing number of contingencies. It can therefore be solved using dense $\L D\L^\transpose$ factorization and back substitution algorithms. The solution of the dense Schur system,
\begin{equation}
	\bfS \Delta \bfu_g =  \b_C - \sum_{i=0}^{N_c} \B_i \A_i^{-1} \b_i, \label{eq:SCsolve}
\end{equation}
yields a part of the solution corresponding to the coupling variables $\Delta \bfu_g$, which is used to obtain all the local solutions $\Delta \bfu_i$ by solving 
\begin{equation}
	\A_i \Delta \bfu_i = \b_i - \B_i ^\transpose\Delta \bfu_g. \label{eq:SCsolve1}
\end{equation}
Since the block contributions to the Schur  $\B_i \A_i^{-1} \B_i^\transpose$ complement are independent,   they can be evaluated in parallel, as well as the residuals $\B_i \A_i^{-1} \b_i$ and the solution $\Delta \bfu_i$ can be computed independently at each process. Interprocess communication occurs because the local Schur complement contributions and Schur complement residuals need to be assembled by the \textit{master} process, and during the broadcast of the Schur complement solution to the remaining  processes. 

In the description of Algorithm \ref{alg:Schur}, sequential steps such as reduction and broadcast are performed only by the \textit{master} process. 


\begin{remark}
One should bear in mind that the computational efficiency obtained by exploiting the block-diagonal structure, such as~\eqref{eq:hessianPerm}, is determined by the number of the coupling variables $|\bfu_g|$. If coupling is large, then the Schur decomposition will not be efficient compared to the direct factorization techniques because of the cubic complexity of dense factorizations \eqref{eq:SCsolve}.
\end{remark}

\begin{algorithm}[t]
	\begin{algorithmic}[1]
		\REQUIRE KKT system with arrowhead structure~\eqref{eq:hessianPerm}, right-hand side $\b$
		\ENSURE $\Delta \bfu$
		\STATE Distribute blocks from the KKT system~\eqref{eq:hessianPerm} evenly across $\mathcal{P}$ processes, where $\mathcal{N}_p$ is the set of diagonal blocks assigned to process $p \in \mathcal{P}$ 
		\STATE Factorize $\A_i = \L_i D_i \L_i^\transpose$ for each $i \in \mathcal{N}_p$\label{alg:SchurSCf1} 
		\STATE Compute $\bfS_i = \B_i \A_i^{-1} \B_i^\transpose$ for each $i \in \mathcal{N}_p$
		\STATE Accumulate $\C_p = \sum_{i \in \mathcal{N}_p} \bfS_i$
		\IF{\textit{master}}
		\STATE Reduce $\bfS = \C - \sum_{p \in \mathcal{P}} \C_p$ \label{alg:SchurSCf2} 
		\ENDIF
		\STATE Compute $\r_i =  \B_i \A_i^{-1} \b_i$ for each $i \in \mathcal{N}_p$ \label{alg:SchurRHS}
		\STATE Accumulate $\r_p =  \sum_{p \in \mathcal{P}} \r_i$ 
		\IF{\textit{master}}
		\STATE Reduce $\r =  \sum_{p \in \mathcal{P}} \r_p$ 
		\STATE Factorize $\bfS = L_s D_s L_s^\transpose$ \label{alg:schurF}
		\STATE Solve $\bfS \Delta \bfu_g =  \b_C - \r$ \label{alg:schurS}
		\STATE Broadcast solution $\bfu_g$ to all $p \in \mathcal{P}$ \label{alg:solve1}
		\ENDIF
		\STATE Solve $\A_i \Delta \bfu_i = \B_i \bfu_g - \b_i $  for each $i \in \mathcal{N}_p$ \label{alg:solve2}
	\end{algorithmic}
	\caption{Parallel procedure for solving the linear systems based on the Schur complement decomposition~\eqref{eq:SCform}--\eqref{eq:SCsolve1}.}
	\label{alg:algo:bw} \label{alg:Schur}
\end{algorithm}

The most expensive step of the presented computational scheme is evaluation of the local contributions to the Schur complement $\B_i \A_i^{-1} \B_i^\transpose$ in~\eqref{eq:SCform}. The standard approach uses a direct sparse solver, such as \PARDISO{}~\cite{pardiso}, to factorize the symmetric matrix $\A_i = \L_iD_i\L_i^\transpose$ and perform multiple forward-backward substitutions with all right-hand side (RHS) vectors in $\B_i^\transpose$, followed by multiplication from the left by $\B_i$. This approach, however, does not exploit sparsity of the problem in $\B_i^\transpose$ blocks, since the linear solver treats the RHS vectors as being dense. 

An alternative approach, implemented in \PARDISO{}~\cite{olafAugmented}, addresses these limitations by performing an incomplete factorization of the augmented matrix $\M_i$:
\begin{equation}
	\M_i = \begin{pmatrix}
	\A_i & \B_i^\transpose \\
	\B_i & 0
	\end{pmatrix},
	\label{eq:augmented} 
\end{equation}
exploiting also the sparsity of $\B_i^\transpose$.  The factorization of $\M_i$ is stopped after pivoting reaches the last diagonal entry of $\A_i$. At this point, the term $-\B_i \A_i^{-1} \B_i^\transpose$ is computed and resides in the (2, 2) block of $\M_i$. By exploiting the sparsity not only in $\A_i$, but also in $\B_i$ it is possible to reduce memory traffic by using in-memory sparse matrix compression techniques, which render this approach quite favorable for multicore parallelization.

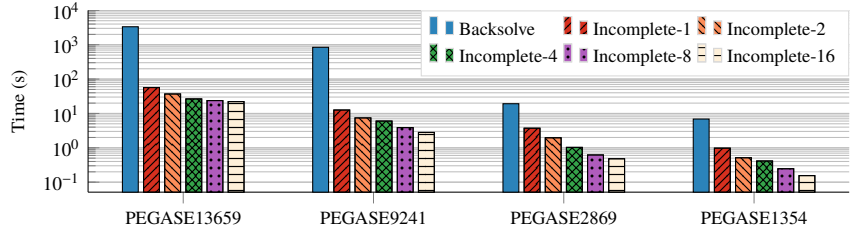
\begin{figure}[t]
    \begin{center}
\begin{tikzpicture}


\pgfplotstableread{
	CASE    Incomplete.01   Incomplete.02   Incomplete.04   Incomplete.08   Incomplete.16   std ID
	PEGASE13659  5.68E+01    3.69E+01    2.65E+01    2.37E+01    2.21E+01    3.34E+03    0
	PEGASE9241   1.26E+01    7.42E+00    6.01E+00    3.88E+00    2.80E+00    8.43E+02    1
	PEGASE2869   3.72E+00    1.96E+00    1.03E+00    0.62E+00    0.48E+00    1.92E+01    2
	PEGASE1354   9.78E-01    5.13E-01    4.14E-01    2.46E-01    1.55E-01    6.88E+00    3
}\datatable;

\begin{axis}[name=symb,
width=\columnwidth, height=4cm, 
ybar,
bar width=6pt, 
ymode=log, 
log origin=infty,
ymax=1e4,
xmin=-0.5,
xmax=3.5,
axis lines*=left, 
ymajorgrids, yminorgrids,
xticklabels from table={\datatable}{CASE},
xtick={0, 1, 2, 3, 4, 5, 6, 7, 8},
xticklabel style={rotate=0, xshift=-0.0cm, anchor=north, font=\scriptsize},
ytick={0.1, 1, 10, 100, 1000, 10000},
yticklabel style={font=\scriptsize},
ylabel style={font=\scriptsize},
xlabel style={font=\scriptsize},
ylabel={Time (s)},
legend style={at={(1.0,1.0)},legend cell align=left,align=right,draw=white!85!black,font=\scriptsize,legend columns=3}
]

\addplot[fill=mycolor0] table[x=ID, y=std] {\datatable};
\addlegendentry{Backsolve}
\addplot[fill=mycolor1, postaction={pattern=north east lines}] table[x=ID, y=Incomplete.01] {\datatable}; 
\addlegendentry{Incomplete-1}
\addplot[fill=mycolor2, postaction={pattern=north west lines}] table[x=ID, y=Incomplete.02] {\datatable}; 
\addlegendentry{Incomplete-2}
\addplot[fill=mycolor3, postaction={pattern=crosshatch}] table[x=ID, y=Incomplete.04] {\datatable}; 
\addlegendentry{Incomplete-4}
\addplot[fill=mycolor4, postaction={pattern=dots}] table[x=ID, y=Incomplete.08] {\datatable}; 
\addlegendentry{Incomplete-8}
\addplot[fill=mycolor5, postaction={pattern=horizontal lines}] table[x=ID, y=Incomplete.16] {\datatable}; 
\addlegendentry{Incomplete-16}
\end{axis}

\end{tikzpicture}
    \end{center}
    \caption{Incomplete factorization of the augmented matrix.}
    \label{fig:augmented}

    \end{figure}

In Figure~\ref{fig:augmented} we compare the standard, so-called backsolve, technique and the multicore incomplete factorization with increasing number of cores is shown for various benchmarks. This demonstrates that the incomplete factorization approach is orders of magnitude faster, especially for the large problems. Due to the extensive memory requirements for storing the RHS vectors in the ``backsolve'' approach, only its single-core execution is demonstrated.

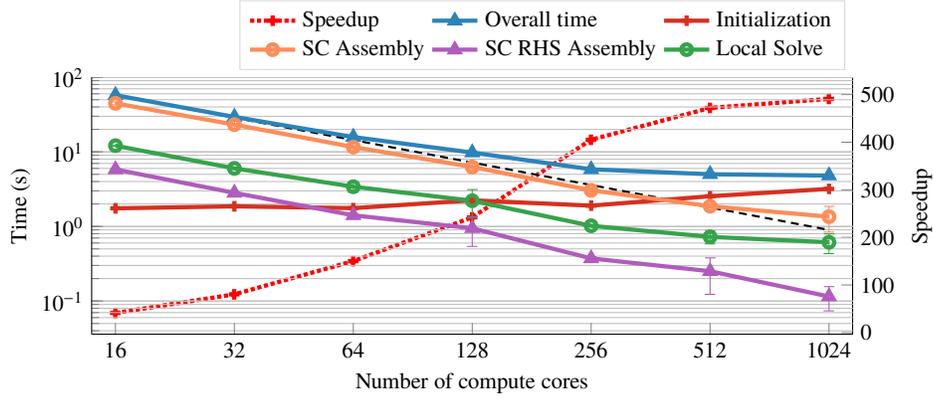
\begin{figure}[!t]
\begin{tikzpicture}

\pgfplotstableread{
ID	NPROC	SSTINITFIRST	SSTINITAVG	SSTSCHURFAVG	SSTSCHURFSTD	SSTSCHURSAVG	SSTSCHURSSTD	SSTSOLVEAVG	SSTSOLVESTD	SSTTOTALAVG	SSTTOTALSTD	IDEAL	SPEEDUP
0	2	2.63E+01	6.69E-01	9.11E+01	1.84E+00	1.14E+01	5.53E-02	1.39E+01	7.42E-02	1.06E+02	5.42E-01	1.06E+02	29.33018868
1	4	2.63E+01	7.66E-01	4.50E+01	8.99E-01	1.01E+01	4.24E-02	1.14E+01	4.86E-02	5.76E+01	7.86E-02	53.2	53.97569444
2	8	2.63E+01	7.61E-01	2.52E+01	5.81E-01	9.47E+00	3.69E-02	1.02E+01	4.02E-02	3.63E+01	4.12E-01	26.6	85.64738292
3	16	2.69E+01	8.25E-01	1.51E+01	4.73E-01	9.27E+00	4.02E-02	9.65E+00	5.21E-02	2.57E+01	4.30E-01	13.3	120.9727626
4	32	2.64E+01	8.65E-01	1.39E+01	1.32E-01	9.03E+00	3.70E-02	9.27E+00	5.18E-02	2.41E+01	1.45E-01	6.65	129.0041494
}\datatablecaseOLDa;

\pgfplotstableread{
ID	NPROC	NTHREADS	SSTINITFIRST	SSTINITAVG	SSTSCHURFAVG	SSTSCHURFSTD	SSTSCHURSAVG	SSTSCHURSSTD	SSTSOLVEAVG	SSTSOLVESTD	SSTTOTALAVG	SSTTOTALSTD	IDEAL	SPEEDUP
0	1	8	4.28E+01	1.01E+00	2.25E+01	5.49E-01	2.85E+00	6.48E-02	5.95E+00	1.14E-01	2.93E+01	5.68E-01	2.93E+01	27.11
1	2	8	4.00E+01	9.27E-01	1.16E+01	1.57E-01	1.47E+00	1.62E-02	3.15E+00	3.46E-02	1.57E+01	1.53E-01	1.47E+01	50.66
2	4	8	4.58E+01	1.48E+00	5.69E+00	2.24E-01	7.39E-01	2.66E-02	1.65E+00	5.88E-02	8.74E+00	1.64E-01	7.33E+00	91.01
3	8	8	4.15E+01	1.44E+00	3.00E+00	4.80E-02	3.78E-01	1.66E-02	8.97E-01	2.92E-02	5.32E+00	1.51E-01	3.67E+00	149.53
4	16	8	4.29E+01	1.22E+00	1.58E+00	5.07E-02	1.95E-01	1.56E-02	5.36E-01	1.83E-02	3.31E+00	3.71E-02	1.83E+00	240.51
5	32	8	4.50E+01	1.50E+00	9.45E-01	6.54E-02	1.00E-01	1.01E-02	3.83E-01	1.26E-02	2.80E+00	1.28E-01	9.16E-01	284.06
6	64	8	4.45E+01	1.90E+00	6.16E-01	3.94E-02	5.22E-02	8.75E-03	2.37E-01	1.05E-02	2.73E+00	4.24E-02	4.58E-01	291.05
}\datatableOLDb;

\pgfplotstableread{
ID	NPROC	NTHREADS	SSTINITFIRST	SSTINITAVG	SSTSCHURFAVG	SSTSCHURFSTD	SSTSCHURSAVG	SSTSCHURSSTD	SSTSOLVEAVG	SSTSOLVESTD	SSTTOTALAVG	SSTTOTALSTD	IDEAL	SPEEDUP
0	1	8	7.92E+01	1.75E+00	4.47E+01	2.06E+00	5.81E+00	1.55E-01	1.21E+01	2.40E-01	5.75E+01	7.15E-01	5.75E+01	4.11E+01
1	2	8	8.65E+01	1.86E+00	2.32E+01	3.12E+00	2.82E+00	8.20E-02	6.01E+00	1.56E-01	2.94E+01	3.41E-01	2.87E+01	8.03E+01
2	4	8	8.88E+01	1.76E+00	1.16E+01	1.40E+00	1.41E+00	3.62E-02	3.41E+00	5.36E-01	1.58E+01	2.38E-01	1.44E+01	1.50E+02
3	8	8	9.33E+01	2.24E+00	6.27E+00	9.59E-01	9.45E-01	4.06E-01	2.21E+00	9.15E-01	9.78E+00	1.70E-01	7.19E+00	2.42E+02
4	16	8	8.74E+01	1.90E+00	3.05E+00	2.29E-01	3.73E-01	1.88E-02	1.02E+00	3.41E-02	5.83E+00	6.55E-02	3.59E+00	4.05E+02
5	32	8	9.33E+01	2.54E+00	1.87E+00	1.20E-01	2.51E-01	1.28E-01	7.25E-01	1.41E-01	5.01E+00	7.06E-02	1.80E+00	4.72E+02
6	64	8	1.04E+02	3.19E+00	1.35E+00	5.09E-01	1.15E-01	4.11E-02	6.13E-01	1.81E-01	4.81E+00	3.13E-02	8.98E-01	4.91E+02
}\datatable;

\begin{axis}[
axis y line*=right,
axis x line=none,
xmin=-0.2,
xmax=6.2,
width=\columnwidth, height=5cm, 
ytick={0,100,200,300,400,500},
yticklabel style={font=\footnotesize},
ylabel style={font=\footnotesize},
ylabel={Speedup},
y label style={at={(axis description cs:1.115,.5)},anchor=south},
y tick label style={
	/pgf/number format/.cd,
	fixed,
	fixed zerofill,
	precision=0,
	/tikz/.cd
},
legend=none
] 
\addplot[color=red, ultra thick, densely dotted, mark=+,  mark options={solid}] table[x=ID, y=SPEEDUP] {\datatable}; \label{plot:speed2}
\end{axis}

\begin{axis}[name=symb,
width=\columnwidth, height=5cm, 
bar width=10pt, 
ymode=log, 
log origin=infty,
ymax=1e2,
xmin=-0.2,
xmax=6.2,
axis lines*=left, 
ymajorgrids, yminorgrids,
xticklabels from table={\datatable}{NPROC},
xtick={0, 1, 2, 3, 4, 5, 6, 7, 8, 9},
xticklabels={16, 32,64,128,256,512,1024,2048,4096,8192,16384},
xticklabel style={rotate=0, xshift=-0.0cm, anchor=north, font=\footnotesize},
yticklabel style={font=\footnotesize},
ylabel style={font=\footnotesize},
y label style={at={(axis description cs:-0.07,.5)},anchor=south},
xlabel style={font=\footnotesize},
ylabel={Time (s)},
xlabel={Number of compute cores},
legend style={at={(0.99,1.3)},legend cell align=left,align=right,draw=white!85!black,font=\footnotesize,legend columns=3}
]
\addlegendimage{/pgfplots/refstyle=plot:speed2}\addlegendentry{Speedup}

\addplot[color=black,  thick, densely dashed, mark options={solid},forget plot] table[x=ID, y=IDEAL] {\datatable};

\addplot[color=mycolor0, ultra thick, mark=triangle, mark options={solid},error bars/.cd, y dir=both, y explicit] table[x=ID, y=SSTTOTALAVG, y error=SSTTOTALSTD] {\datatable}; 
\addlegendentry{Overall time}

\addplot[color=mycolor1, ultra thick, mark=+, mark options={solid},error bars/.cd, y dir=both, y explicit] table[x=ID, y=SSTINITAVG] {\datatable}; 
\addlegendentry{Initialization}
\addplot[color=mycolor2, ultra thick, mark=o, mark options={solid},error bars/.cd,    y dir=both, y explicit] table[x=ID, y=SSTSCHURFAVG, y error=SSTSCHURFSTD] {\datatable}; 
\addlegendentry{SC Assembly}
\addplot[color=mycolor4, ultra thick, mark=triangle, mark options={solid},error bars/.cd,    y dir=both, y explicit] table[x=ID, y=SSTSCHURSAVG, y error=SSTSCHURSSTD] {\datatable}; 
\addlegendentry{SC RHS Assembly}
\addplot[color=mycolor3, ultra thick, mark=o, mark options={solid},error bars/.cd, y dir=both, y explicit] table[x=ID, y=SSTSOLVEAVG, y error=SSTSOLVESTD] {\datatable}; 
\addlegendentry{Local Solve}  
\end{axis}

\end{tikzpicture}
\caption{Scaling of the parallel approach using the PEGASE1354-4096 benchmark and the speedup with respect to the direct sparse solver.}
\label{fig:scaling}
\end{figure}

We evaluated the strong scaling efficiency of the distributed solver on the ``Piz Daint" supercomputer, using an increasing number of compute cores on the distributed compute nodes. The instance of the solved problem contained up to $1.1\cdot10^7$ variables and $2.7\cdot10^7$ constraints and the size of the KKT system is $5.48\cdot10^7$.
Figure \ref{fig:scaling} shows the average wall time of the individual phases of Algorithm~\ref{alg:Schur}, indicating also the ideal strong scaling of the overall time. The algorithmic phases presented are the initialization phase, assembly of the Schur complement using the incomplete factorization of the augmented matrix in steps~\ref{alg:SchurSCf1} --~\ref{alg:SchurSCf2}, RHS vector assembly and Schur complement solution in steps~\ref{alg:schurF} --~\ref{alg:schurS}, and solutions of the local parts of the system in steps~\ref{alg:solve1} --~\ref{alg:solve2}. Figure~\ref{fig:scaling} also demonstrates the speedups of the distributed solver compared to the serial direct factorization. The benchmarks were run with a single MPI process per node and 16 threads per process.

The distributed approach using a single process outperforms the sequential direct factorization by a factor of up to $40\times$. With an increasing number of distributed nodes the observed speedup was up to $500\times$.
The distributed solution time scales reasonably up to $512$ cores at $32$ compute nodes, which in terms of workload translates to $128$ scenarios per node of PEGASE1354 benchmark. At this point, the most expensive part of the algorithm, the computation of the local contributions to the Schur complement, requires approximately the same time as the initialization phase, where the KKT system is distributed to all available compute nodes.
The acceleration and efficiency of the structure-exploiting algorithm stems from the reduced complexity associated with the factorization of the smaller sparse diagonal blocks instead of the original  SCOPF KKT system \eqref{eq:SCOPFhessianSymmetric} or its reduced variant \eqref{eq:SCOPFhessianReduced} after the it was permuted to the arrowhead structure \eqref{eq:hessianPerm}.
For sufficiently large power grids, however, the dense Schur complement (SC) system might become very large, and dominate the overall processing time in steps \ref{alg:schurF} and \ref{alg:schurS}. Hardware accelerators such as GPUSs might be deployed to address the computational complexity of the dense linear algebra. Otherwise, the dimensions of the dense systems remain feasible for the majority of power grids, since the dimensions depend only on the power grid properties, not on the number of contingency scenarios.

\subsection{Structure-Exploiting Algorithms  for MPOPF}
 For the MPOPF problems, the size of the dense SC grows very quickly, not only with the size of the network but also proportionally to the number of installed storage devices and the number of time periods $N\Ns$.  As the number of time periods $N$ or storage devices $\Ns$ increases, the solution approach based on Algorithm \ref{alg:Schur} results in a less efficient algorithm than the  direct sparse approach employing PARDISO on the original KKT system \eqref{eq:SCOPFhessian}, both with respect to computational time and memory consumption despite the benefits of the Schur decomposition. However,  the MPOPF problem, unlike the SCOPF problem, can be optimized even further by exploiting the particular structure of the off-diagonal blocks $\B_n$.

Inspecting the particular structure of the blocks $\B_n$ \eqref{Eq:B12N}, one can see that the  SC matrix computed by \eref{eq:augmented} for the $n$th  block $\Sc_n = - \B_n \A_n^{-1} \B_n^\transpose$  has the structure
\begin{align}
	\label{Eq:LocalSchurBs}
	\Sc_n = \pmat{  \O_n         &\0^\transpose_n       &\0^\transpose_n                     & \cdots& \0^\transpose_n                  \\
		\0_n         &{\color{mycolor0} \Sc_{11,n}}     & {\color{mycolor0} \Sc_{10,n}^\transpose}       & \cdots& \Sc_{10,n}^\transpose            \\ 
		\0_n         &{\color{mycolor0} \Sc_{10,n}}            & {\color{mycolor0} \Sc_{00,n}}                  & \cdots& \Sc_{00,n}^\transpose            \\
		\vdots       &\vdots              & \vdots                & \ddots& \vdots                \\
		\0_n         &\Sc_{10,n}            & \Sc_{00,n}              & \cdots& \Sc_{00,n}},
\end{align}
where the $\0_n \in \R{\Ns \times (n-1) \Ns }$, $\O_n \in \R{(n-1) \Ns \times (n-1) \Ns }$, and $\Sc_{ij,n} = - \C_i \A_n^{-1} \C_j^\transpose$, $\; i,j \in \{0, 1\}$.
The only blocks in $\Sc_n$ that are distinct are colored in  blue and form the entries of the 
2 by 2 block matrix
\begin{align}
	\label{Eq:LocalSchurBs3}
	\discr{\bar{S}}_n = \pmat{\Sc_{11,n}&            \Sc_{10,n}^\transpose \\ 
		\Sc_{10,n}&            \Sc_{00,n}},
\end{align}
where the rest of the rows and columns of $\Sc_n$ are direct replicates of the entries of 
the last row and column of $\discr{\bar{S}}_n$.

\begin{figure}[!t]
\begin{subfigure}[t]{0.5\textwidth}
\centering
\newcommand\offone{-8.0ex}
\newcommand\offtwo{-10.5ex}
\newcommand\offthree{-4.5ex}

\begin{tikzpicture}[scale=1.0]

\begin{axis}
[
    separate axis lines,
    ymajorgrids,
    yminorgrids,
    minor tick num=1,
    thick,
    every axis plot post/.style={/pgf/number format/fixed},
    ybar=2pt,
    bar width=4pt,
    xmin=3600,
    ymin=1,
    ymax=1e6,
    xtick=data,
    ymode=log,
    log basis y = 10,
    width=1.05\columnwidth,
    yminorticks=true,
    enlarge x limits=0.17,
    separate axis lines,
    every outer x axis line/.append style={black},
    every x tick label/.append style={font=\scriptsize},
    every y tick label/.append style={font=\scriptsize},
    symbolic x coords={3600,4800,7200,8760},
    restrict y to domain*=1:1e6, 
    visualization depends on=rawy\as\rawy, 
    after end axis/.code={ 
            \draw [ultra thick, white, decoration={snake, amplitude=1pt}, decorate] (rel axis cs:0,1.01) -- (rel axis cs:1,1.01);
        },
    axis lines*=left,
    clip=false,
    axis background/.style={fill=none},
    legend style={at={(0.95,1.22)},legend cell align=left,align=left,draw=none,fill=none,font=\scriptsize, legend columns=2},
    xlabel/.append style={font=\scriptsize},
    ylabel/.append style={font=\scriptsize},
    xlabel={$N$},
    ylabel={Time (s)}
]
\addplot[fill=mycolor1]                                          table [x=N, y=AveragePerIter]{data/figures/HYPEROPTGeneMem0Fig4.dat}; \addlegendentry{BELTISTOS}
\addplot[fill=mycolor2, postaction={pattern=horizontal lines}]   table [x=N, y=AveragePerIter]{data/figures/HYPEROPTGeneMem1Fig4.dat}; \addlegendentry{BELTISTOSmem}      
\addplot[fill=mycolor3, postaction={pattern=crosshatch}]         table [x=N, y=AveragePerIter]{data/figures/IPOPTGeneFig4.dat}; \addlegendentry{IPOPT}
\addplot[fill=mycolor4, postaction={pattern=north east lines}]   table [x=N, y=AveragePerIter]{data/figures/MIPSGeneFig4.dat}; \addlegendentry{MIPS}
\addplot[fill=mycolor5, postaction={pattern=north west lines}]   table [x=N, y=AveragePerIter]{data/figures/KNITROGeneFig4.dat}; \addlegendentry{KNITRO}
\end{axis}
\end{tikzpicture}
\caption{Average time per iteration.}
\label{fig:AllSolversAveragePerformance}
\end{subfigure}%
~
\begin{subfigure}[t]{0.5\textwidth}
\centering
\begin{tikzpicture}[scale=1.0]

\begin{axis}
[
    separate axis lines,
    ymajorgrids,
    yminorgrids,
    minor tick num=0,
    thick,
    every axis plot post/.style={/pgf/number format/fixed},
    ybar=3pt,
    bar width=4pt,
    width=1.05\columnwidth,
    xmin=3600,
    ymin=1,
    ymax=1e6,
    ymode=log,
    log basis y = 10,
    yminorticks=true,
    enlarge x limits=0.17,
    every outer x axis line/.append style={black},
    every x tick label/.append style={font=\scriptsize},
    every y tick label/.append style={font=\scriptsize},
    symbolic x coords={3600,4800,7200,8760},
    restrict y to domain*=1:1e6, 
    visualization depends on=rawy\as\rawy, 
    after end axis/.code={ 
            \draw [ultra thick, white, decoration={snake, amplitude=1pt}, decorate] (rel axis cs:0,1.05) -- (rel axis cs:1,1.05);
        },
    restrict y to domain*=1:1e6, 
    axis lines*=left,
    clip=false,
    axis background/.style={fill=none},
    legend style={at={(0.95,1.22)},legend cell align=left,align=left,draw=none,fill=none, font=\scriptsize,legend columns=2},
    xlabel/.append style={font=\scriptsize},
    ylabel/.append style={font=\scriptsize},
    xlabel={$N$},
    ylabel={Memory (MB)}
]
\addplot[fill=mycolor1]                                       table [x=N, y=MEMORY]{data/figures/HYPEROPTGeneMem0Fig4.dat};  \addlegendentry{BELTISTOS}
\addplot[fill=mycolor2,postaction={pattern=horizontal lines}] table [x=N, y=MEMORY]{data/figures/HYPEROPTGeneMem1Fig4.dat};  \addlegendentry{BELTISTOSmem}
\addplot[fill=mycolor3,postaction={pattern=crosshatch} ]      table [x=N, y=MEMORY]{data/figures/IPOPTGeneFig4.dat};     \addlegendentry{IPOPT}

\addlegendimage{empty legend}
\addlegendentry{}

\addplot[fill=mycolor5,postaction={pattern=north west lines}] table [x=N, y=MEMORY]{data/figures/KNITROGeneFig5.dat}; 
\addlegendentry{KNITRO}
\end{axis}
\end{tikzpicture}
\caption{Memory allocated.}
\label{fig:MemoryChart}
\end{subfigure}%
\caption{Case IEEE118. Statistics for solving the KKT system.}
\label{fig:MPOPF}
\end{figure}
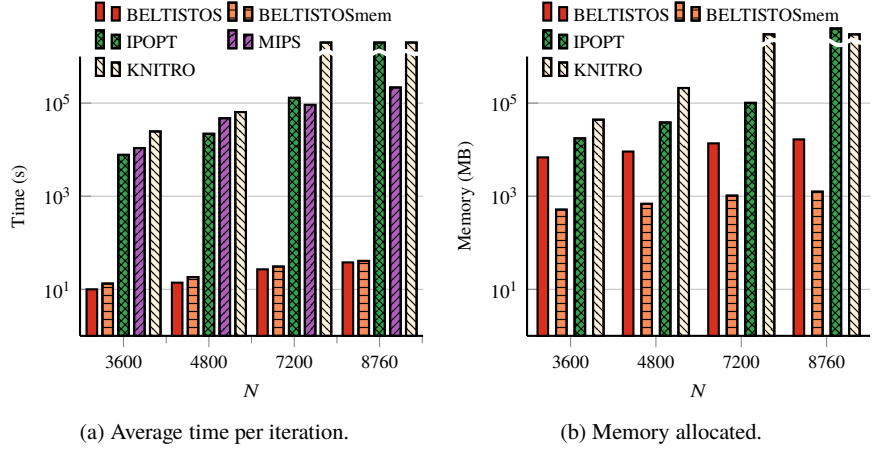

Since each one of the blocks of $\discr{\bar{S}}_n$ has size $\Ns \times \Ns$, the computation of $\Sc_n$ becomes independent of the number of time periods $N$
and only depends on the number of storage devices $\Ns$. It is easily verified that  the global SC $\Sc_c$ obtains the form
\begin{align}
	\label{Eq:GlobalSchur}
	\Sc_c = \pmat{\Sc_{11}   & \Sc_{12}^\transpose&  \Sc_{12}^\transpose   &\cdots  & \Sc_{12}^\transpose \\ 
		\Sc_{12}   & \Sc_{22}           &  \Sc_{23}^\transpose   &\cdots  & \Sc_{23}^\transpose \\ 
		\Sc_{12}   & \Sc_{23}           &  \Sc_{33}              &\cdots  & \Sc_{34}^\transpose \\
		\vdots     & \vdots             & \vdots                 &\ddots  & \vdots   \\
		\Sc_{12}   & \Sc_{23}           & \Sc_{34}               &\cdots  & \Sc_{NN}},
\end{align}
where each block of $\Sc_c \in \R{N \Ns \times N \Ns}$ has dimensions $\Ns \times \Ns$. 
Storing $\Sc_c$ due to its special structure, requires only two block vectors:
one for all diagonal blocks $\Sc_d = [\Sc_{11},\Sc_{22}, \cdots, \Sc_{NN}]$ of size $\Ns \times N \Ns$,
and one for the off-diagonal blocks $\Sc_o = [\Sc_{12},\Sc_{23}, \cdots, \Sc_{N-1N}]$ of size $\Ns \times (N-1) \Ns$,
significantly reducing this way the storage requirements for $\Sc_c$. Furthermore, exploiting the fact
that the blocks below the main diagonal of each column of $\Sc_c$ in
\eref{Eq:GlobalSchur} are identical, we can perform the factorization in
$O(n^2)$ operations instead of $O(n^3)$, which is the case for standard dense $\L D\L^\transpose$ factorization of $\Sc_c$ with $n=N \Ns$. Similarly, the back substitution can be performed in $O(n)$ instead of $O(n^2)$.  The reduction in the computational
complexity and storage requirements of the SC system
renders the overall approach significantly more 
economical in terms of overall running time and memory footprint, as demonstrated in Figure \ref{fig:MPOPF}.

For comparison,  we also consider three alternative
optimization algorithms that also adopt an IP strategy, namely,
\IPOPT{}~\cite{IPOPT2005,IPOPT2006},
\MIPS{}~\cite{Wang:2007}, and \KNITRO~\cite{KNITRO}. The structure exploiting IP algorithm introduced in this section is referred to as \BELTISTOS{}.
The average time per iteration for $N=3600$ up to $N=8760$ corresponding to one
year with a time step size corresponding to one hour, is shown in Figure
\ref{fig:AllSolversAveragePerformance}.  For this set of benchmarks \KNITRO{}
needed more than 1 TB of memory for $N \ge 5760$ and it terminated with a related
error message. For $N=8760$ \PARDISO{} failed due to overflow of the
number of nonzero entries in the $\L,D$ factors. It is worth noting that
\BELTISTOS{mem} (the memory saving approach of \BELTISTOS{} that implements
Algorithm~\ref{alg:Schur} without storing the factors of the blocks $\A_i$ and computing them on the fly
in steps \ref{alg:SchurSCf1}, \ref{alg:SchurRHS}, and \ref{alg:solve2}), although it is slightly slower than the
normal mode of \BELTISTOS{}, it is still almost four orders of magnitude faster
than \IPOPT{} and \MIPS{}. It also needs approximately two orders of magnitude less
memory than \IPOPT{} as it is shown in figure \ref{fig:MemoryChart}, where we plot
the memory (in MB) allocated by each algorithm for the solution of the KKT system.
The \MIPS{} and \KNITRO{} solvers do not report the memory allocated and it could only be estimated for the case of \KNITRO{}.

\section{Results and Discussion}

This study demonstrates that significant performance gains are possible, for specific classes of optimal control problems,
not by exploiting supercomputers and parallel distributed or multithreaded programming, but through deeper understanding of the problem structure and the design of algorithms adapted to the problem structure. Orders of magnitude of faster execution time and orders of magnitude of memory savings were achieved rendering the solution of very-large-scale 
problems, previously intractable without a supercomputer, possible on a common laptop~\cite{beltistos}.

The Schur decomposition enables low memory SC assembly on a per-block basis, whenever a problem can be reordered to an arrowhead structure, which is the case for many real life problems composed of enumerated subproblems, such as contingency scenarios for SCOPF problems or time periods for MPOPF problems, while at the same time promoting parallel processing.
%
Even on single-core execution for SCOPF problems, speedups from 40 fold to 270 fold were observed while further exploitation of distributed multicore and many-core computing environments for the solution of the structured KKT system drastically reduces the execution times and demonstrates significant progress towards the solution of large-scale SCOPF problems.

In contrast to SCOPF and although MPOPF problems can be reordered into an arrowhead structure, the reordering results in a  dense SC that grows in size with the number of time periods and does not necessarily lead to a more efficient solution strategy. However, owing to the intrinsic structure of the linear constraints, the Schur decomposition     algorithm supplemented with elimination strategies exploiting data compression, resulted in an overall solution strategy of unprecedented performance. Memory was reduced by approximately two orders of magnitude, while runtime performance still remains about three orders of magnitude higher than competitors, even on a single core. 

Our findings strongly motivate further structural inspection and analysis of the present and  similar problems of the same family, anticipating that adopting and extending the presented structure-exploiting techniques for other problems would result in significant acceleration of other OPF problems of interest paving the way for the next generation of OPF algorithms.

\bibliographystyle{abbrv}
\begin{small}
\bibliography{biblio}
\end{small}
\end{document}